\documentclass[12pt]{iopart}
\usepackage{amssymb}
\usepackage{iopams}
\usepackage{graphicx}
\usepackage{xcolor}
\usepackage{footmisc}

\definecolor{myblue}{rgb}{0.000,0.000,1.000}

\eqnobysec
\begin{document}

\title[Convexification With Viscocity for EIT]{Convexification With the
Viscocity Term for Electrical Impedance Tomography}
\author{Michael V. Klibanov$^{1}$, Jingzhi Li$^{2}$ and Zhipeng Yang$^{3}$}

\begin{abstract}
A version of the globally convergent convexification numerical method is
constructed for the problem of Electrical Impedance Tomography in the 2D
case. An important element of this version is the presence of the viscosity
term. Global convergence analysis is carried out. Results of numerical
experiments are presented.
\end{abstract}

%

\address{$^1$ Department of Mathematics and Statistics, University of North
Carolina at Charlotte, Charlotte, NC, 28223, USA} 
\address{$^2$ Department of Mathematics, Southern University of Science
and Technology, Shenzhen 518055, P.~R.~China} 
\address{$^3$ Department of
School of Mathematics and Statistics, Lanzhou University, Lanzhou 730000, P.
R. China} 
\eads{\mailto{mklibanv@charlotte.edu},
\mailto{li.jz@sustech.edu.cn}, \mailto{yangzp@lzu.edu.cn}}

\vspace{10pt} \begin{indented}
\item[]
\end{indented}




\maketitle

%
%

\section{Introduction}

\label{sec:1}

The electrical impedance tomography (EIT) problem has gained a significant
popularity in the Inverse Problems community. We cite now some related
publications. Since this paper is not a survey, then we restrict citations
only to a few references \cite%
{Gehre,Hamilton,Harrach1,Harrach2,Harrach3,KEIT,KL,Mueller}. The EIT is a
non expensive imaging method aimed to recover the spatially distributed
electrical conductivity inside of an object of interest using boundary
measurements. EIT has direct applications in medical imaging \cite%
{Holder,Mueller}.

\textbf{Definition. }\emph{We call a numerical method for a Coefficient
Inverse Problem globally convergent if there exists a theorem which ensures
that this method delivers points in a sufficiently small neighborhood of the
true solution of this problem without any advanced knowledge of this
neighborhood. If, however, convergence is guaranteed only if the starting
point of iterations is located in that neighborhood, then we call this
method locally convergent one.}

This paper is concerned with the development of a version of the globally
convergent convexification numerical method for a 2D Coefficient Inverse
Problem (CIP) for EIT. More precisely, we work with such a version of the
convexification method for the considered CIP, which uses a viscosity term.
This was not done for our CIP in the past. The paper \cite{KlibHJ} is the
originating publication, in which the idea of the introduction of the
viscosity term in the convexification method was presented for the first
time. The goal of \cite{KlibHJ} was to solve numerically the Hamilton-Jacobi
equation. In addition, we refer to the more recent publication \cite{LN} on
this subject. As to the use of the viscosity term in the convexification
method for CIPs, the idea of \cite{KlibHJ,LN} was extended in \cite{Ktransp}
to a CIP for the radiative transport equation.

The CIP of this paper is about the reconstruction of the unknown spatially
dependent electric conductivity coefficient from boundary measurements. In 
\cite{KEIT} and \cite[Chapter 7]{KL} a version of the convexification method
without the viscosity term was developed for our CIP. The presence of the
viscosity term allows us to obtain a boundary value problem for a system of
only two coupled elliptic PDEs to solve. On the other hand, in \cite{KEIT,KL}
a system of $N>2$ coupled nonlinear elliptic equations was used. That system
was generated by an expansion of the solution of a certain PDE derived from
the original one in a truncated Fourier-like series with respect to a
special orthonormal basis. That basis was originally introduced in \cite%
{KJIIP}. The number $N$ was the number of terms in that truncated series.

The work \cite{KI} is the first one, where the convexification method was
introduced in the field of Coefficient Inverse Problems with the goal to
avoid the well known phenomenon of multiple local minima and ravines of
conventional least squares mismatch functionals for these problems, see,
e.g. \cite{B1,B2,B3,Chavent,Giorgi,Gonch1,Gonch2,Grote,LB,Riz} and
references cited therein for those functionals. The convexification method
constructs globally convergent numerical procedures for CIPs, unlike locally
convergent ones of these references. On the other hand, the global
convergence issue for EIT is discussed in publications \cite%
{Harrach1,Harrach2,Harrach3}. More precisely, the paper \cite{Harrach1}
presents a numerical example of multiple local minima for a CIP\ for EIT. To
focus on the global convergence issue, \cite{Harrach1} contains a quite
fresh idea allowing to equivalently reformulate the EIT problem with a
finite number of measurements as a convex semidefinite problem. In addition,
it is shown in \cite{Harrach2} how to globalize the level-set method for EIT
by starting it with the globally convergent monotonicity method. The work 
\cite{Harrach3} is about the closely related Robin problem, where, in
addition to the global convergence, a rigorous characterization is provided
of the number of electrodes, which are required to ensure a certain
resolution.

While originally the convexification was proposed in \cite{KI} for a CIP
only for a hyperbolic PDE, currently a number of versions of this method are
developed for CIPs for PDEs for hyperbolic, parabolic, elliptic, transport
PDEs, as well as for the system of PDEs governing mean field games, see,
e.g. \cite{KJIIP,SAR,KEIT,KL,Ktr,Ktransp,MFG} and references cited therein.
Each version of the convexification requires its own convergence analysis
and its own numerical implementation. Convergence rates are explicitly
written for all versions.

The convexification consists of two steps. On the first step, the original
CIP is transformed to a boundary value problem (BVP) for such a nonlinear
PDE (or a system of coupled PDEs), which does not contain the unknown
coefficient. The boundary conditions are Dirichlet and Neumann boundary
conditions on at least a part of the boundary of the domain of interest. On
the second step, a globally strongly convex weighted Tikhonov-like
functional is constructed to solve that BVP. The weight is the Carleman
Weight Function (CWF), i.e. the function, which is involved in the Carleman
estimate for the corresponding PDE operator. The global convergence of the
gradient descent method of the minimization of that functional is
established and explicit formulae for convergence rates of iterations of
that method to the true solution of the CIP are given. Convergence to the
true solution takes place as long as the noise in the input data tends to
zero.

In section 2 we formulate both forward and inverse problems. In section 3 we
present our transformation procedure. In section 4 we present the
convexification method. In section 5 theorems of the convergence analysis
are formulated. These theorems are proven in sections 6 and 7. In section 8
we present results of our numerical experiments.

\section{Statements of Forward and Inverse Problems}

\label{sec:2}

Denote $\mathbf{x}=\left( x,y\right) $ points in $\mathbb{R}^{2}.$ Let $%
G\subset \mathbb{R}^{2}$ be a bounded domain with a sufficiently smooth
boundary $\partial G.$ Let $\left( a,b\right) \subset G$ be a point in $G$
and let $A>0$ be a number. Let $D_{A}\left( a,b\right) $ be the disk of the
radius $A$ with the center at $\left( a,b\right) $ and let $C_{A}\left(
a,b\right) $ be the circle, which is the boundary of this disk,%
\begin{equation}
\left. 
\begin{array}{c}
D_{A}\left( a,b\right) =\left\{ \mathbf{x}=\left( x,y\right) :\left(
x-a\right) ^{2}+\left( y-b\right) ^{2}<A^{2}\right\} , \\ 
C_{A}\left( a,b\right) =\left\{ \mathbf{x}=\left( x,y\right) :\left(
x-a\right) ^{2}+\left( y-b\right) ^{2}=A^{2}\right\} .%
\end{array}%
\right.  \label{2.1}
\end{equation}%
We assume that 
\begin{equation}
D_{A}\left( a,b\right) \subset G,\mbox{ }C_{A}\left( a,b\right) \cap
\partial G=\varnothing .  \label{2.2}
\end{equation}%
Let $c_{1},c_{2}>0$ be two numbers such that 
\begin{equation}
\left. 
\begin{array}{c}
\Omega =\left\{ \mathbf{x}=\left( x,y\right) :\left\vert x-a\right\vert
<c_{1},\left\vert y-b\right\vert <c_{2}\right\} \subset D_{A}\left(
a,b\right) , \\ 
\left\{ x=0\right\} \cap \overline{\Omega }=\varnothing .%
\end{array}%
\right.  \label{2.20}
\end{equation}%
be a rectangle with the center at the point $\left( a,b\right) .$ Then $%
\partial \Omega \cap C_{A}\left( a,b\right) =\varnothing .$ Let $\Gamma
\subset \partial \Omega $ be the side of the rectangle $\Omega $ defined as:%
\begin{equation}
\Gamma =\left\{ \mathbf{x}=\left( x,y\right) :x=a+c_{1},\mbox{ }y\in \left(
b-c_{2},b+c_{2}\right) \right\}  \label{2.21}
\end{equation}%
We assume that point sources $\mathbf{x}_{0}\mathbf{\in }C_{A}\left(
a,b\right) $ run along the circle $C_{A}\left( a,b\right) .$ It is natural
to model the point source as the delta function $\delta \left( \mathbf{x}-%
\mathbf{x}_{0}\right) .$ However, to simplify the presentation, we model the
point source as the $\delta -$like function $f\left( \mathbf{x}-\mathbf{x}%
_{0}\right) $ \cite{KEIT}, \cite[page 141, formula (7.2)]{KL}, 
\begin{equation}
f\left( \mathbf{x-x}_{0}\right) =C_{\rho }\left\{ 
\begin{array}{c}
\exp \left( \frac{\left\vert \mathbf{x-x}_{0}\right\vert ^{2}}{\rho
^{2}-\left\vert \mathbf{x-x}_{0}\right\vert ^{2}}\right) ,\left\vert \mathbf{%
\ x-x}_{0}\right\vert <\rho , \\ 
0,\left\vert \mathbf{x-x}_{0}\right\vert \geq \rho%
\end{array}%
\right.  \label{2.3}
\end{equation}%
where $\rho \in \left( 0,1\right) $ is a certain number. The number $C_{\rho
}$ is chosen such that%
\[
C_{\rho }\int_{\left\vert \mathbf{x}\right\vert <\rho }\exp \left( \frac{%
\left\vert \mathbf{x}\right\vert ^{2}}{\rho ^{2}-\left\vert \mathbf{x}%
\right\vert ^{2}}\right) d\mathbf{x}=1, 
\]%
and also $\rho $ is so small that 
\begin{equation}
f\left( \mathbf{x}-\mathbf{x}_{0}\right) =0\mbox{ for }\mathbf{x}\in 
\overline{\Omega }.  \label{2.4}
\end{equation}

Let $u\left( \mathbf{x},\mathbf{x}_{0}\right) $ be the electrical potential
at the point $\mathbf{x}\in G$ generated by the point source $\mathbf{x}%
_{0}\in C_{A}\left( a,b\right) .$ Let $\sigma \left( \mathbf{x}\right) $ be
the electrical conductivity of the medium at the point $\mathbf{x.}$ We
assume that%
\begin{equation}
\sigma \in C^{2}\left( \overline{G}\right) ,  \label{2.5}
\end{equation}%
\begin{equation}
\sigma \left( \mathbf{x}\right) \geq 1\mbox{ in }G\mbox{,}  \label{2.6}
\end{equation}%
\begin{equation}
\sigma \left( \mathbf{x}\right) =1\mbox{ in }G\diagdown \Omega .  \label{2.7}
\end{equation}%
The equation of the EIT is%
\begin{equation}
\mbox{div}\left( \sigma \left( \mathbf{x}\right) \nabla u\right) =-f\left( 
\mathbf{x}-\mathbf{x}_{0}\right) ,\mbox{ }\mathbf{x}\in G.  \label{2.8}
\end{equation}%
We impose the zero Dirichlet boundary condition at $\partial G,$%
\begin{equation}
u\left( \mathbf{x},\mathbf{x}_{0}\right) \mid _{\mathbf{x}\in \partial G}=0.
\label{2.9}
\end{equation}%
The forward problem of EIT is the problem of finding the solution $u\left( 
\mathbf{x},\mathbf{x}_{0}\right) \in C^{2}\left( \overline{G}\right) $ for $%
\mathbf{x}_{0}\in C_{A}\left( a,b\right) $ of the Dirichlet boundary value
problem (BVP) (\ref{2.8}), (\ref{2.9}). It follows from the maximum
principle and (\ref{2.1})-(\ref{2.9}) that there exists a number $d>0$ such
that \cite[Chapter 3, Theorem 3.1]{GT}%
\begin{equation}
u\left( \mathbf{x},\mathbf{x}_{0}\right) \geq d,\mbox{ }\forall \left( 
\mathbf{x},\mathbf{x}_{0}\right) \in \Omega \times C_{A}\left( a,b\right) .
\label{2.10}
\end{equation}%
Our focus here is on the numerical solution of the following CIP:

\textbf{Coefficient Inverse Problem (CIP).} \emph{Let conditions (\ref{2.1}
)-(\ref{2.9})} be in place. \emph{Suppose that functions }$g_{0}$\emph{\ and 
}$g_{1}$\emph{\ are known, where}%
\begin{equation}
\left. 
\begin{array}{c}
u\left( \mathbf{x},\mathbf{x}_{0}\right) \mid _{\mathbf{x}\in \partial
\Omega }=g_{0}\left( \mathbf{x},\mathbf{x}_{0}\right) ,\mbox{ }\forall
\left( \mathbf{x},\mathbf{x}_{0}\right) \in \partial \Omega \times
C_{A}\left( a,b\right) , \\ 
\partial _{n}u\left( \mathbf{x},\mathbf{x}_{0}\right) \mid _{\mathbf{x}\in
\Gamma }=g_{1}\left( \mathbf{x},\mathbf{x}_{0}\right) ,\mbox{ }\forall 
\mathbf{x}\in \Gamma ,\mbox{ }\forall \mathbf{x}_{0}\in C_{A}\left(
a,b\right) ,%
\end{array}
\right.  \label{2.11}
\end{equation}%
\emph{where }$\Gamma \subset \partial \Omega $\emph{\ is defined in (\ref%
{2.21})} \emph{and }$\partial _{n}$\emph{\ is the normal derivative with
respect to }$\mathbf{x}$\emph{. Find the coefficient }$\sigma \left( \mathbf{%
\ x}\right) $\emph{.}

Hence, it is assumed in this CIP that the Dirichlet boundary condition is
known on the entire boundary $\partial \Omega $ of the rectangle $\Omega $
and, in addition, the Neumann boundary condition is known on the side $%
\Gamma $ of this rectangle.

\section{Transformation}

\label{sec:3}

Change variables via introducing the new function $p\left( \mathbf{x},%
\mathbf{x}_{0}\right) ,$ 
\begin{equation}
p\left( \mathbf{x},\mathbf{x}_{0}\right) =\sqrt{\sigma \left( \mathbf{x}
\right) }u\left( \mathbf{x},\mathbf{x}_{0}\right) .  \label{3.1}
\end{equation}%
Then by (\ref{2.3})-(\ref{2.8}), (\ref{2.11}) and (\ref{3.1}) 
\begin{equation}
\Delta p+a\left( \mathbf{x}\right) p=0,\mbox{ }\mathbf{x}\in \Omega ,
\label{3.2}
\end{equation}%
\begin{equation}
\left. 
\begin{array}{c}
p\left( \mathbf{x},\mathbf{x}_{0}\right) \mid _{\mathbf{x}\in \partial
\Omega }=g_{0}\left( \mathbf{x},\mathbf{x}_{0}\right) ,\mbox{ }\forall
\left( \mathbf{x},\mathbf{x}_{0}\right) \in \partial \Omega \times
C_{A}\left( a,b\right) , \\ 
\partial _{n}p\left( \mathbf{x},\mathbf{x}_{0}\right) \mid _{\mathbf{x}\in
\Gamma }=g_{1}\left( \mathbf{x},\mathbf{x}_{0}\right) ,\mbox{ }\forall 
\mathbf{x}\in \Gamma ,\forall \mathbf{x}_{0}\in C_{A}\left( a,b\right) .%
\end{array}
\right.  \label{3.3}
\end{equation}%
\begin{equation}
a\left( \mathbf{x}\right) =\frac{\Delta \left( \sqrt{\sigma \left( \mathbf{x}
\right) }\right) }{\sqrt{\sigma \left( \mathbf{x}\right) }}.  \label{3.4}
\end{equation}

Using (\ref{2.10}) and (\ref{3.1}), consider the second change of variables
as%
\begin{equation}
v\left( \mathbf{x},\mathbf{x}_{0}\right) =\ln \left( p\left( \mathbf{x},%
\mathbf{x}_{0}\right) \right) .  \label{3.5}
\end{equation}%
Using (\ref{3.2}), (\ref{3.3}) and (\ref{3.5}), we obtain%
\begin{equation}
\Delta v+\left( \nabla v\right) ^{2}+a\left( \mathbf{x}\right) =0,\mbox{ }%
\mathbf{x}\in \Omega ,  \label{3.6}
\end{equation}%
\begin{equation}
v\left( \mathbf{x},\mathbf{x}_{0}\right) \mid _{\mathbf{x}\in \partial
\Omega }=s_{0}\left( \mathbf{x},\mathbf{x}_{0}\right) ,\mbox{ }\forall
\left( \mathbf{x},\mathbf{x}_{0}\right) \in \partial \Omega \times
C_{A}\left( a,b\right) ,  \label{3.7}
\end{equation}%
\begin{equation}
\partial _{n}v\left( \mathbf{x},\mathbf{x}_{0}\right) \mid _{\mathbf{x}\in
\Gamma }=s_{1}\left( \mathbf{x},\mathbf{x}_{0}\right) ,\mbox{ }\forall 
\mathbf{x}\in \Gamma ,\forall \mathbf{x}_{0}\in C_{A}\left( a,b\right) ,
\label{3.8}
\end{equation}%
where 
\begin{equation}
s_{0}\left( \mathbf{x},\mathbf{x}_{0}\right) =\ln \left( g_{0}\left( \mathbf{%
\ x},\mathbf{x}_{0}\right) \right) ,  \label{3.9}
\end{equation}%
\begin{equation}
s_{1}\left( \mathbf{x},\mathbf{x}_{0}\right) =\frac{g_{1}\left( \mathbf{x},%
\mathbf{x}_{0}\right) }{g_{0}\left( \mathbf{x},\mathbf{x}_{0}\right) }.
\label{3.10}
\end{equation}%
Using polar coordinates, we obtain for $\mathbf{x}_{0}\in C_{A}\left(
a,b\right) $%
\begin{equation}
\left. 
\begin{array}{c}
\mathbf{x}_{0}=\left( x_{0},y_{0}\right) ,\mbox{ }x_{0}=a+A\cos \varphi ,%
\mbox{ }y_{0}=b+A\sin \varphi , \\ 
\varphi \in \left( 0,2\pi \right) .%
\end{array}%
\right.   \label{3.11}
\end{equation}%
Hence, we replace below the notation $\mathbf{x}_{0}$ with $\varphi .$ Set 
\begin{equation}
w\left( \mathbf{x},\varphi \right) =v\left( \mathbf{x},\mathbf{x}_{0}\right)
,\mbox{ }\mathbf{x}\in \Omega ,\varphi \in \left( 0,2\pi \right) .
\label{3.12}
\end{equation}%
To eliminate the unknown coefficient $a\left( \mathbf{x}\right) $ from
equation (\ref{3.6}) we differentiate this equation with respect to $\varphi
.$ We obtain%
\begin{equation}
\Delta w_{\varphi }+2\nabla w_{\varphi }\nabla w=0,\mathbf{x}\in \Omega
,\varphi \in \left( 0,2\pi \right) ,  \label{3.13}
\end{equation}%
\begin{equation}
w_{\varphi }\mid _{\mathbf{x}\in \partial \Omega }=\partial _{\varphi
}s_{0}\left( \mathbf{x},\varphi \right) ,\mbox{ }\forall \left( \mathbf{x}%
,\varphi \right) \in \partial \Omega \times \left( 0,2\pi \right) ,
\label{3.14}
\end{equation}%
\begin{equation}
\partial _{\varphi }\left( \partial _{n}w\left( \mathbf{x},\varphi \right)
\right) \mid _{\mathbf{x}\in \Gamma }=\partial _{\varphi }s_{1}\left( 
\mathbf{x},\varphi \right) ,\mbox{ }\forall \left( \mathbf{x},\varphi
\right) \in \Gamma \times \left( 0,2\pi \right) ,  \label{3.15}
\end{equation}%
where functions $s_{0}$ and $s_{1}$ are those defined in (\ref{3.9}), (\ref%
{3.10}).

\ An inconvenience of equation (\ref{3.13}) is that it contains two unknown
functions $w$ and $w_{\varphi }.$ In the case of the convexification method
for a CIP for time dependent equations, such as ones in, e.g. \cite{KI,MFG}, 
\cite[chapter 9]{KL}, $w_{\varphi }$ would be replaced with $w_{t}$ and then 
$w$ would be expressed via $w_{t}$ using an initial condition for $w$ and a
Volterra integral with respect to $t$. In our case, however, we do not have
an initial condition for $w.$

Hence, a truncated Fourier-like series with respect to a special orthonormal
basis of functions depending on $\varphi $ was used in \cite{KEIT}, \cite[%
chapter 7]{KL} to approximately solve BVP (\ref{3.13})-(\ref{3.15}).\ This
basis was first introduced in \cite{KJIIP}. Now, however, we introduce in
equation (\ref{3.13}) the viscosity term $-\varepsilon \Delta w$ instead of
using that basis. More precisely we replace equation (\ref{3.13}) with 
\begin{equation}
-\varepsilon \Delta w+\Delta w_{\varphi }+2\nabla w_{\varphi }\nabla w=0,%
\mathbf{x}\in \Omega ,\mbox{ }\varphi \in \left( 0,2\pi \right) ,
\label{3.16}
\end{equation}%
where $\varepsilon \in \left( 0,1\right) $ is a sufficiently small number to
be chosen computationally. The problem of a rigorous proof of the
convergence of our process as $\varepsilon \rightarrow 0$ is a substantially
challenging one. Therefore, we do not address this problem here.

Denote%
\begin{equation}
q=w_{\varphi }-\varepsilon w.  \label{3.17}
\end{equation}%
Then 
\begin{equation}
w=\frac{w_{\varphi }-q}{\varepsilon }.  \label{3.18}
\end{equation}%
Then (\ref{3.13}) and (\ref{3.16}) become:%
\begin{equation}
L_{1}\left( w_{\varphi },q\right) =\Delta w_{\varphi }+2\nabla w_{\varphi
}\nabla \left( \frac{w_{\varphi }-q}{\varepsilon }\right) =0,\mathbf{x}\in
\Omega ,\varphi \in \left( 0,2\pi \right) ,  \label{3.19}
\end{equation}%
\begin{equation}
L_{2}\left( w_{\varphi },q\right) =\Delta q+2\nabla w_{\varphi }\nabla
\left( \frac{w_{\varphi }-q}{\varepsilon }\right) =0,\mathbf{x}\in \Omega
,\varphi \in \left( 0,2\pi \right) .  \label{3.20}
\end{equation}%
Boundary conditions for the system of PDEs (\ref{3.19}), (\ref{3.20}) can be
derived from (\ref{3.7})-(\ref{3.10}), (\ref{3.12}), (\ref{3.14}), (\ref%
{3.15}), (\ref{3.17}) and (\ref{3.18}):%
\begin{equation}
w_{\varphi }\left( \mathbf{x},\varphi \right) \mid _{\mathbf{x}\in \partial
\Omega }=\partial _{\varphi }s_{0}\left( \mathbf{x},\varphi \right) ,\mbox{ }%
\left( \mathbf{x},\varphi \right) \in \partial \Omega \times \left( 0,2\pi
\right) ,  \label{3.21}
\end{equation}%
\begin{equation}
\partial _{n}\left( \partial _{\varphi }w\left( \mathbf{x},\varphi \right)
\right) \mid _{\mathbf{x}\in \Gamma }=\partial _{\varphi }s_{1}\left( 
\mathbf{x},\varphi \right) ,\mbox{ }\forall \left( \mathbf{x},\varphi
\right) \in \Gamma \times \left( 0,2\pi \right) ,  \label{3.22}
\end{equation}%
\begin{equation}
q\left( \mathbf{x},\varphi \right) \mid _{\mathbf{x}\in \partial \Omega
}=\partial _{\varphi }s_{0}\left( \mathbf{x},\varphi \right) -\varepsilon
s_{0}\left( \mathbf{x},\varphi \right) ,\forall \left( \mathbf{x},\varphi
\right) \in \partial \Omega \times \left( 0,2\pi \right) ,  \label{3.23}
\end{equation}%
\begin{equation}
\partial _{n}q\left( \mathbf{x},\varphi \right) =\partial _{\varphi
}s_{1}\left( \mathbf{x},\varphi \right) -\varepsilon s_{1}\left( \mathbf{x}%
,\varphi \right) ,\mbox{ }\forall \left( \mathbf{x},\varphi \right) \in
\Gamma \times \left( 0,2\pi \right) .  \label{3.24}
\end{equation}

Thus, we focus below on the solution of BVP (\ref{3.19})-(\ref{3.24}) with
respect to the pair $\left( w,q\right) $. Our solution will depend on $%
\varphi $ as a parameter. Suppose that we have solved this BVP. Then we will
calculate our approximate unknown coefficient using (\ref{3.6}), (\ref{3.12}%
) and (\ref{3.18}) as%
\begin{equation}
a\left( \mathbf{x}\right) =-\frac{1}{2\pi }\int\limits_{0}^{2\pi }\left(
\Delta w+\left( \nabla w\right) ^{2}\right) \left( \mathbf{x},\varphi
\right) d\varphi ,\mbox{ }\mathbf{x}\in \Omega .  \label{3.25}
\end{equation}

\section{Convexification}

\label{sec:4}

Introduce spaces $H$ and $H_{2},$ 
\begin{equation}
\left. 
\begin{array}{c}
H=\left\{ 
\begin{array}{c}
p\left( \mathbf{x},\varphi \right) :p\left( \mathbf{x},\varphi \right) \in
H^{3}\left( \Omega \right) ,\forall \varphi \in \left[ 0,2\pi \right] , \\ 
\left\Vert p\right\Vert _{H}=\sup_{\varphi \in \left( 0,2\pi \right) }\left(
\left\Vert p\left( \mathbf{x},\varphi \right) \right\Vert _{H^{3}\left(
\Omega \right) }\right) <\infty ,%
\end{array}%
\right\} \\ 
H_{2}=\left\{ \left( r,s\right) \left( \mathbf{x},\varphi \right) \in
H\times H\right\} , \\ 
\left\Vert \left( r,s\right) \left( \mathbf{x},\varphi \right) \right\Vert
_{H_{2}}^{2}=\left\Vert r\left( \mathbf{x},\varphi \right) \right\Vert
_{H}^{2}+\left\Vert s\left( \mathbf{x},\varphi \right) \right\Vert _{H}^{2}.%
\end{array}%
\right.  \label{4.1}
\end{equation}%
Let $R>0$ be an arbitrary number. Define the set $B\left( R\right) \subset
H_{2}$ 
\begin{equation}
B\left( R\right) =\left\{ 
\begin{array}{c}
\left( r,s\right) \left( \mathbf{x},\varphi \right) \in H_{2}: \\ 
\left\Vert \left( r,s\right) \right\Vert _{H_{2}}<R,\mbox{ } \\ 
r\mbox{ satisfies boundary conditions (\ref{3.21}), (\ref{3.22}),} \\ 
s\mbox{ satisfies boundary conditions (\ref{3.23}), (\ref{3.24}).}%
\end{array}%
\right\}  \label{4.2}
\end{equation}%
By the Sobolev embedding theorem and (\ref{4.2}) 
\begin{equation}
\hspace{-1 cm}
\left. 
\begin{array}{c}
r\left( \mathbf{x},\varphi \right) ,s\left( \mathbf{x},\varphi \right) \in
C^{1}\left( \overline{\Omega }\right) ,\forall \varphi \in \left( 0,2\pi
\right) ,\mbox{ }\forall \left( r,s\right) \in \overline{B\left( R\right) },
\\ 
\sup_{\varphi \in \left( 0,2\pi \right) }\left( \left\Vert r\left( \mathbf{x}%
,\varphi \right) \right\Vert _{C^{1}\left( \overline{\Omega }\right)
}^{2}+\left\Vert s\left( \mathbf{x},\varphi \right) \right\Vert
_{C^{1}\left( \overline{\Omega }\right) }^{2}\right) \leq C,\mbox{ }\forall
\left( r,s\right) \in \overline{B\left( R\right) },%
\end{array}%
\right.  \label{4.02}
\end{equation}%
where $C=C\left( R\right) >0$ is a certain number depending only on $R$.

The Carleman Weight Function is the same as the one in a slightly modified
formula (10.29) of \cite{KL}:%
\begin{equation}
\psi _{\lambda }\left( \mathbf{x}\right) =\exp \left( 2\lambda x^{2}\right) .
\label{4.3}
\end{equation}%
Introduce the functional $J_{\lambda ,\alpha }\left( r,s\right) \left(
\varphi \right) ,$ which depends on $\varphi $ as on a parameter,%
\[
J_{\lambda ,\alpha }:\overline{B\left( R\right) }\rightarrow \mathbb{R}, 
\]%
\begin{equation}
\hspace{-1cm}\left. 
\begin{array}{c}
J_{\lambda ,\alpha }\left( r,s\right) \left( \varphi \right) =\sqrt{%
\varepsilon }\int\limits_{\Omega }\left[ \left( L_{1}\left( r,s\right)
\left( \mathbf{x},\varphi \right) \right) ^{2}+\left( L_{2}\left( r,s\right)
\left( \mathbf{x},\varphi \right) \right) ^{2}\right] \psi _{\lambda }\left( 
\mathbf{x}\right) d\mathbf{x}+ \\ 
+\alpha \left\Vert \left( r,s\right) \left( \mathbf{x},\varphi \right)
\right\Vert _{H^{3}\left( \Omega \right) \times H^{3}\left( \Omega \right)
}^{2},\mbox{ }\forall \varphi \in \left( 0,2\pi \right) ,%
\end{array}%
\right.  \label{4.4}
\end{equation}%
where the nonlinear operators $L_{1}$ and $L_{2}$ are defined in (\ref{3.19}%
) and (\ref{3.20}) respectively and $\alpha \in \left( 0,1\right) $ is the
regularization parameter. We introduce the multiplier $\sqrt{\varepsilon }$
in (\ref{4.4}) for two reasons. First, to balance two terms in the right
hand side of (\ref{4.4}) since by (\ref{2.20}) and (\ref{4.3}) $\max_{%
\overline{\Omega }}\psi _{\lambda }\left( \mathbf{x}\right) =\exp \left[
2\lambda \left( a+c_{1}\right) ^{2}\right] ,$ whereas $\alpha \in \left(
0,1\right) .$ Second, we have computationally established that the term $%
\sqrt{\varepsilon }$ in (\ref{4.4}) works well for computations, see section
8. We consider the following problem:

\textbf{Minimization Problem 1.} \emph{Minimize the functional }$J_{\lambda
}\left( r,s\right) \left( \varphi \right) $\emph{\ on the set }$\overline{
B\left( R\right) }$ \emph{for each} \emph{value of the parameter} $\varphi
\in \left( 0,2\pi \right) .$

Suppose that the minimizer $\left( r_{\min },s_{\min }\right) \left( \mathbf{%
\ x},\varphi \right) \in \overline{B\left( R\right) }$ of functional (\ref%
{4.4}) is found. Then we set by (\ref{3.18}) and (\ref{3.19})%
\begin{equation}
\left. 
\begin{array}{c}
w_{\varphi }\left( \mathbf{x},\varphi \right) =r_{\min }\left( \mathbf{x}
,\varphi \right) ,\mbox{ }w\left( \mathbf{x},\varphi \right) =\left(
1/\varepsilon \right) \left( r_{\min }-s_{\min }\right) \left( \mathbf{x}
,\varphi \right) , \\ 
\left( \mathbf{x},\varphi \right) \in \Omega \times \left( 0,2\pi \right) .%
\end{array}
\right.  \label{4.5}
\end{equation}%
Next, the function $a\left( \mathbf{x}\right) $ is found via (\ref{3.25}).

As soon as the function $a\left( \mathbf{x}\right) $ is found, it seems to
be that the function $\sqrt{\sigma \left( \mathbf{x}\right) }$ can be found
as the solution of the following Dirichlet boundary value problem, which
follows from (\ref{2.5})-(\ref{2.7}) and (\ref{3.4}):%
\begin{equation}
\Delta \left( \sqrt{\sigma \left( \mathbf{x}\right) }\right) -a\left( 
\mathbf{x}\right) \left( \sqrt{\sigma \left( \mathbf{x}\right) }\right) =0,
\label{4.6}
\end{equation}%
\begin{equation}
\sqrt{\sigma \left( \mathbf{x}\right) }\mid _{\partial \Omega }=1.
\label{4.7}
\end{equation}%
However, the existence and uniqueness theorem for this problem require $%
a\left( \mathbf{x}\right) \geq 0$ \cite{GT,LU}. This inequality is not
guaranteed in the above procedure. Therefore, we use the Neumann boundary
condition as well, which follows from (\ref{2.5})-(\ref{2.7}),%
\begin{equation}
\partial _{n}\left( \sqrt{\sigma \left( \mathbf{x}\right) }\right) \mid _{ 
\mathbf{x}\in \partial \Omega }=0.  \label{4.8}
\end{equation}

BVP (\ref{4.6})-(\ref{4.8}) has the overdetermination in the boundary
conditions. Hence, we find an approximate solution of BVP (\ref{4.6})-(\ref%
{4.8}) via the variational Quasi-Reversibility Method (QRM) \cite[chapter 4]%
{KL}, which is basically a linear version of the Minimization Problem 1.
More precisely, we solve the Minimization Problem 2.

\textbf{Minimization Problem 2.} \emph{Minimize the quadratic functional }$%
K\left( V\right) ,$%
\begin{equation}
K\left( V\right) =\int\limits_{\Omega }\left( \Delta V-a\left( \mathbf{x}%
\right) V\right) ^{2}d\mathbf{x,}  \label{4.9}
\end{equation}%
\emph{over the set of functions }$V\in H^{2}\left( \Omega \right) $\emph{\
satisfying the boundary conditions:}%
\begin{equation}
V\mid _{\mathbf{x}\in \partial \Omega }=1, \quad \partial _{n}V\mid _{%
\mathbf{x}\in \partial \Omega }=0.  \label{4.10}
\end{equation}

Suppose that we have found the minimizer $V_{\min }\left( \mathbf{x}\right) $
of functional (\ref{4.9}) with the boundary conditions (\ref{4.10}). Then we
set%
\begin{equation}
\sigma \left( \mathbf{x}\right) =V_{\min }^{2}\left( \mathbf{x}\right) .
\label{4.11}
\end{equation}%
Existence and uniqueness of the minimizer $V_{\min }\in H^{2}\left( \Omega
\right) $ satisfying boundary conditions (\ref{4.10}) can be proven
similarly with Theorem 4.3.1 of \cite{KL}. Hence, a convergence analysis for
the Minimization Problem 2 is omitted here for brevity.

\section{Theorems of the Convergence Analysis}

\label{sec:5}

In this section we formulate four theorems of our convergence analysis. One
of them, about a Carleman estimate is known, another one we prove in this
section, and two others are proven in section 6.

\subsection{The global strong convexity}

\label{sec:5.1}

We start from a Carleman estimate. Denote 
\begin{equation}
\left. 
\begin{array}{c}
H_{0}^{2}\left( \Omega \right) =\left\{ u\in H^{2}\left( \Omega \right)
:u\mid _{\partial \Omega }=0,\partial _{n}u\mid _{\Gamma }=0\right\} , \\ 
H_{0}^{3}\left( \Omega \right) =\left\{ u\in H^{3}\left( \Omega \right)
:u\mid _{\partial \Omega }=0,\partial _{n}u\mid _{\Gamma }=0\right\} .%
\end{array}
\right.  \label{5.1}
\end{equation}%
Obviously $H_{0}^{3}\left( \Omega \right) \subset H_{0}^{2}\left( \Omega
\right) .$

\textbf{Theorem 5.1} (Carleman estimate \cite[Theorem 10.3.1]{KL}). \emph{\
Let }$\psi _{\lambda }\left( \mathbf{x}\right) $\emph{\ be the CWF defined
in (\ref{4.3}). Then there exists a sufficiently large number }$\lambda
_{0}=\lambda _{0}\left( \Omega \right) \geq 1$\emph{\ and a number }$%
C=C\left( \Omega \right) ,$\emph{\ both numbers depending only on }$\Omega ,$%
\emph{\ such that the following Carleman estimate holds:}%
\[
\int\limits_{\Omega }\left( \Delta u\right) ^{2}\psi _{\lambda }\left( 
\mathbf{x}\right) d\mathbf{x\geq }\frac{C}{\lambda }\int\limits_{\Omega
}\left( u_{xx}^{2}+u_{xy}^{2}+u_{yy}^{2}\right) \psi _{\lambda }\left( 
\mathbf{x}\right) d\mathbf{x+} 
\]%
\[
+C\lambda \int\limits_{\Omega }\left( \nabla u\right) ^{2}\psi _{\lambda
}\left( \mathbf{x}\right) d\mathbf{x+}C\lambda ^{3}\int\limits_{\Omega
}u^{2}\psi _{\lambda }\left( \mathbf{x}\right) d\mathbf{x,} 
\]%
\[
\forall \lambda \geq \lambda _{0},\mbox{ }\forall u\in H_{0}^{2}\left(
\Omega \right) . 
\]

The condition in the second line of (\ref{2.20}) ensures that $\left\vert
\nabla \psi _{\lambda }\left( \mathbf{x}\right) \right\vert \neq 0$ in $%
\overline{\Omega },$ which is important for the proof of this theorem.

\textbf{Theorem 5.2} (the central result). \emph{The following holds true:}

\emph{1. For each }$\lambda >0$\emph{, for each }$\varphi \in \left( 0,2\pi
\right) $\emph{\ and for each pair }$\left( r,s\right) \in \overline{B\left(
R\right) }$\emph{\ the functional }$J_{\lambda ,\alpha }\left( r,s\right)
\left( \varphi \right) $\emph{\ has the Fr\'{e}chet derivative }%
\begin{equation}
J_{\lambda ,\alpha }^{\prime }\left( r,s\right) \left( \varphi \right) \in
H_{0}^{3}\left( \Omega \right) \times H_{0}^{3}\left( \Omega \right) \emph{,}
\forall \varphi \in \left( 0,2\pi \right) .  \label{5.01}
\end{equation}%
\emph{\ Furthermore, this derivative is Lipschitz continuous on }$\overline{
B\left( R\right) },$\emph{\ i.e. there exists a number }$D=D\left( \lambda
,\varphi ,R,\Omega \right) >0$\emph{\ depending only on listed parameters
such that} 
\begin{equation}
\left. 
\begin{array}{c}
\left\Vert J_{\lambda ,\alpha }^{\prime }\left( r_{2},s_{2}\right) \left(
\varphi \right) -J_{\lambda ,\alpha }^{\prime }\left( r_{1},s_{1}\right)
\left( \varphi \right) \right\Vert _{H^{3}\left( \Omega \right) \times
H^{3}\left( \Omega \right) }\leq \\ 
\leq D\left\Vert \left( r_{2},s_{2}\right) \left( \mathbf{x},\varphi \right)
-\left( r_{1},s_{1}\right) \left( \mathbf{x},\varphi \right) \right\Vert
_{H^{3}\left( \Omega \right) \times H^{3}\left( \Omega \right) }, \\ 
\mbox{ }\forall \left( r_{1},s_{1}\right) \left( \mathbf{x},\varphi \right)
,\left( r_{2},s_{2}\right) \left( \mathbf{x},\varphi \right) \in \overline{
B\left( R\right) }.%
\end{array}
\right.  \label{5.2}
\end{equation}

\emph{2. Let }$\lambda _{0}=\lambda _{0}\left( \Omega \right) \geq 1$\emph{\
be the number of Theorem 5.1. There exists a sufficiently large number }$%
\lambda _{1}=\lambda _{1}\left( R,\Omega ,\varepsilon \right) \geq \lambda
_{0}$\emph{\ such that for each }$\lambda \geq \lambda _{1}$\emph{\ the
functional }$J_{\lambda ,\alpha }\left( r,s\right) $\emph{\ is strongly
convex on the set }$\overline{B\left( R\right) },$\emph{\ i.e. there exists
a number }$C_{1}=C_{1}\left( R,\Omega ,\varepsilon \right) >0$\emph{\ such
that }%
\begin{equation}
\hspace{-2 cm} \left. 
\begin{array}{c}
J_{\lambda ,\alpha }\left( r_{2},s_{2}\right) \left( \varphi \right)
-J_{\lambda ,\alpha }\left( r_{1},s_{1}\right) \left( \varphi \right)
-J_{\lambda ,\alpha }^{\prime }\left( r_{1},s_{1}\right) \left( \varphi
\right) \left( r_{2}-r_{1},s_{2}-s_{1}\right) \left( \mathbf{x},\varphi
\right) \geq \\ 
\geq C_{1}\exp \left( \lambda \left( a-c_{1}\right) \right) \left\Vert
\left( r_{2}-r_{1},s_{2}-s_{1}\right) \left( \mathbf{x},\varphi \right)
\right\Vert _{H^{2}\left( \Omega \right) \times H^{2}\left( \Omega \right)
}^{2}+ \\ 
+\alpha \left\Vert \left( r_{2}-r_{1},s_{2}-s_{1}\right) \left( \mathbf{x}
,\varphi \right) \right\Vert _{H^{3}\left( \Omega \right) \times H^{3}\left(
\Omega \right) }^{2}, \\ 
\forall \varphi \in \left( 0,2\pi \right) ,\mbox{ }\forall \left(
r_{1},s_{1}\right) \left( \mathbf{x},\varphi \right) ,\left(
r_{2},s_{2}\right) \left( \mathbf{x},\varphi \right) \in \overline{B\left(
R\right) },\mbox{ }\forall \lambda \geq \lambda _{1}.%
\end{array}
\right.  \label{5.3}
\end{equation}%
\emph{Both numbers }$\lambda _{1}\left( R,\Omega ,\varepsilon \right) $\emph{%
\ \ and }$C_{1}\left( R,\Omega ,\varepsilon \right) $\emph{\ depend only on }%
$R,\Omega $\emph{\ and }$\varepsilon .$

\emph{3. For each }$\lambda \geq \lambda _{1}$\emph{\ there exists unique
minimizer }$\left( r_{\min ,\lambda },s_{\min ,\lambda }\right) \left( 
\mathbf{x},\varphi \right) \in \overline{B\left( R\right) }$\emph{\ of the
functional }$J_{\lambda }\left( r,s\right) \left( \varphi \right) $\emph{\
on the set }$\overline{B\left( R\right) }.$\emph{\ Furthermore, the
following inequality holds:}%
\begin{equation}
\hspace{-1 cm} \left. 
\begin{array}{c}
J_{\lambda ,\alpha }^{\prime }\left( r_{\min ,\lambda ,\alpha },s_{\min
,\lambda ,\alpha }\right) \left( \varphi \right) \left( \left( r_{\min
,\lambda ,\alpha },s_{\min ,\lambda ,\alpha }\right) \left( \mathbf{x}
,\varphi \right) -\left( r,s\right) \left( \mathbf{x},\varphi \right)
\right) \leq 0,\mbox{ } \\ 
\forall \varphi \in \left( 0,2\pi \right) ,\forall \left( r,s\right) \left( 
\mathbf{x},\varphi \right) \in \overline{B\left( R\right) }.%
\end{array}
\right.  \label{5.4}
\end{equation}

Below $C_{1}=C_{1}\left( R,\Omega ,\varepsilon \right) >0$ denotes different
numbers depending only on listed parameters.

\subsection{The accuracy of the minimizer}

\label{sec:5.2}

We now present the concept of this subsection. Assume that Theorem 5.1 is
proven. In particular, this theorem guarantees the existence and uniqueness
of the minimizer of functional (\ref{4.4}) on the set \emph{\ }$\overline{%
B\left( R\right) }.$ We now want to estimate the accuracy of this minimizer.
To do this, we follow one of main concepts of the theory of Ill-Posed
Problems \cite{TA,T} by assuming the existence of the exact solution $\sigma
^{\ast }\left( \mathbf{x}\right) $ of our CIP satisfying conditions (\ref%
{2.5})-(\ref{2.7}) and with the \textquotedblleft ideal", i.e. noiseless
data $g_{0}^{\ast }\left( \mathbf{x,}\varphi \right) ,g_{1}^{\ast }\left( 
\mathbf{x,}\varphi \right) $ in (\ref{2.11}). The function $\sigma ^{\ast
}\left( \mathbf{x}\right) $ generates the coefficient $a^{\ast }\left( 
\mathbf{x}\right) $ by (\ref{3.4}) as well as the exact pair of functions $%
\left( r^{\ast },s^{\ast }\right) \left( \mathbf{x},\varphi \right) .$
However, the input data for inverse problems always contain noise. Hence, we
need to assume the existence of the noisy data $g_{0}\left( \mathbf{x,}%
\varphi \right) ,g_{1}\left( \mathbf{x,}\varphi \right) $ in (\ref{2.11})
with a small level of noise $\xi \in \left( 0,1\right) .$ Then we need to
estimate the distance between the minimizer $\left( r_{\min ,\lambda
},s_{\min ,\lambda }\right) \left( \mathbf{x},\varphi \right) $ of the
functional $J_{\lambda }\left( r,s\right) \left( \mathbf{x},\varphi \right) $
and the pair of functions $\left( r^{\ast },s^{\ast }\right) \left( \mathbf{x%
},\varphi \right) $. Using this estimate and (\ref{3.25}), we will estimate
the distance between the exact function $a^{\ast }\left( \mathbf{x}\right) $
and the function $a\left( \mathbf{x}\right) ,$ which we will obtain via $%
\left( r_{\min ,\lambda },s_{\min ,\lambda }\right) \left( \mathbf{x}%
,\varphi \right) .$ If taking into account Minimization Problem 2, then a
further effort would lead to an estimate of the distance between $\sigma
^{\ast }\left( \mathbf{x}\right) $ and the function $\sigma \left( \mathbf{x}%
\right) $ obtained from $a\left( \mathbf{x}\right) $ due to (\ref{4.9})-(\ref%
{4.11}). However, we are not doing the latter for brevity since we do not
provide a theory for Minimization Problem 2, see the last paragraph of
section 4.

Assume that there exists a vector function $P\left( \mathbf{x,}\varphi
\right) ,$ 
\[
P\left( \mathbf{x,}\varphi \right) =\left( P_{1}\left( \mathbf{x,}\varphi
\right) ,P_{2}\left( \mathbf{x,}\varphi \right) \right) \in H_{2} 
\]%
such that $P_{1}\left( \mathbf{x,}\varphi \right) $ satisfies boundary
conditions (\ref{3.21}), (\ref{3.22}) and $P_{2}\left( \mathbf{x,}\varphi
\right) $ satisfies boundary conditions (\ref{3.23}), (\ref{3.24}). Here $%
P_{1}$ stands for $w_{\varphi },$ and $P_{2}$ stands for $q$. Let 
\[
P^{\ast }\left( \mathbf{x,}\varphi \right) =\left( P_{1}^{\ast },P_{2}^{\ast
}\right) \left( \mathbf{x,}\varphi \right) \in H_{2} 
\]%
be a vector function satisfying boundary conditions (\ref{3.21})-(\ref{3.24}%
), in which the vector function $\left( \partial _{\varphi }s_{0},\partial
_{\varphi }s_{1}\right) \left( \mathbf{x},\varphi \right) $ is replaced with 
$\left( \partial _{\varphi }s_{0}^{\ast },\partial _{\varphi }s_{1}^{\ast
}\right) \left( \mathbf{x},\varphi \right) $ and such that 
\begin{equation}
\left\Vert \left( P_{1},P_{2}\right) \left( \mathbf{x,}\varphi \right)
-\left( P_{1}^{\ast },P_{2}^{\ast }\right) \left( \mathbf{x,}\varphi \right)
\right\Vert _{H_{2}}<\xi ,  \label{5.5}
\end{equation}%
where $\xi \in \left( 0,1\right) $ is a sufficiently small number
representing the noise level in the data. Here $\left( \partial _{\varphi
}s_{0}^{\ast },\partial _{\varphi }s_{1}^{\ast }\right) \left( \mathbf{x}%
,\varphi \right) $ is the vector function $\left( \partial _{\varphi
}s_{0},\partial _{\varphi }s_{1}\right) \left( \mathbf{x},\varphi \right) ,$
which corresponds to the exact solution with noiseless data. In addition, we
assume that 
\begin{equation}
\left\Vert P\right\Vert _{H_{2}},\left\Vert P^{\ast }\right\Vert _{H_{2}}<R.
\label{5.6}
\end{equation}%
Denote 
\begin{equation}
\overline{B^{\ast }\left( R\right) }=\left\{ 
\begin{array}{c}
Q\left( \mathbf{x,}\varphi \right) =\left( Q_{1},Q_{2}\right) \left( \mathbf{%
\ x,}\varphi \right) \in H_{2}:\left\Vert Q\right\Vert _{H_{2}}\leq R, \\ 
\mbox{functions }Q_{1}\mbox{ and }Q_{2}\mbox{ satisfy boundary conditions}
\\ 
\mbox{ (\ref{3.21}), (\ref{3.22}) and } \\ 
\mbox{(\ref{3.23}),(\ref{3.24}) respectively,} \\ 
\mbox{in which the pair }\left( s_{0},s_{1}\right) \left( \mathbf{x,}\varphi
\right) \\ 
\mbox{ is replaced with \ the pair }\left( s_{0}^{\ast },s_{1}^{\ast
}\right) \left( \mathbf{x,}\varphi \right) .%
\end{array}
\right\} .  \label{5.7}
\end{equation}%
Hence, the exact pair of functions 
\begin{equation}
\left( r^{\ast },s^{\ast }\right) \left( \mathbf{x},\varphi \right) \in 
\overline{B^{\ast }\left( R\right) }.  \label{5.8}
\end{equation}

For every vector function $\left( r,s\right) \in \overline{B\left( R\right) }
$ consider the vector function 
\begin{equation}
\left( \widetilde{r},\widetilde{s}\right) =\left( r-P_{1},s-P_{2}\right) .
\label{5.9}
\end{equation}%
In addition, introduce the vector function $\left( \widetilde{r}^{\ast },%
\widetilde{s}^{\ast }\right) ,$ 
\begin{equation}
\left( \widetilde{r}^{\ast },\widetilde{s}^{\ast }\right) =\left( r^{\ast
}-P_{1}^{\ast },s^{\ast }-P_{2}^{\ast }\right) .  \label{5.10}
\end{equation}%
Using (\ref{5.6})-(\ref{5.10}) and triangle inequality we obtain 
\begin{equation}
\left( \widetilde{r},\widetilde{s}\right) ,\left( \widetilde{r}^{\ast }, 
\widetilde{s}^{\ast }\right) \in B_{0}\left( 2R\right) ,  \label{5.100}
\end{equation}%
where%
\begin{equation}
\hspace{-1 cm} \overline{B_{0}\left( 2R\right) }=\left\{ 
\begin{array}{c}
Q\left( \mathbf{x,}\varphi \right) =\left( Q_{1},Q_{2}\right) \left( \mathbf{%
\ x,}\varphi \right) \in H_{2}:\left\Vert Q\right\Vert _{H_{2}}\leq 2R, \\ 
\mbox{functions }Q_{1}\mbox{ and }Q_{2}%
\mbox{ satisfy zero boundary
			conditions} \\ 
\mbox{ (\ref{3.21}), (\ref{3.22})} \\ 
\mbox{and (\ref{3.23}), (\ref{3.24}) respectively.}%
\end{array}
\right\}  \label{5.11}
\end{equation}%
Introduce the functional $F_{\lambda }\left( r,s\right) $ 
\begin{equation}
\left. 
\begin{array}{c}
F_{\lambda ,\alpha }\left( r,s\right) :\overline{B_{0}\left( 2R\right) }
\rightarrow \mathbb{R}, \\ 
F_{\lambda ,\alpha }\left( r,s\right) =J_{\lambda ,\alpha }\left(
r+P_{1},s+P_{2}\right) .%
\end{array}
\right.  \label{5.12}
\end{equation}%
It follows from (\ref{5.12}) that an obvious analog of Theorem 5.1 holds for
the functional $F_{\lambda }\left( r,s\right) .$ The triangle inequality, (%
\ref{5.6}) and (\ref{5.100}) imply that 
\[
\left( r+P_{1},s+P_{2}\right) \in \overline{B\left( 3R\right) },\mbox{ }
\forall \left( r,s\right) :\overline{B_{0}\left( 2R\right) }. 
\]%
Hence, in that analog of Theorem 5.1 for the functional $F_{\lambda ,\alpha
}\left( r,s\right) ,$ we should replace $\lambda _{1}\left( R,\Omega
,\varepsilon \right) $ with 
\begin{equation}
\lambda _{2}=\lambda _{1}\left( 3R,\Omega ,\varepsilon \right) .
\label{5.13}
\end{equation}

\textbf{Theorem 5.2} (the accuracy of the minimizer). \emph{Assume that
conditions (\ref{5.5})-(\ref{5.13}) hold. For any }$\lambda \geq \lambda
_{2},$\emph{\ let }$\left( \widetilde{r}_{\min ,\lambda },\widetilde{s}%
_{\min ,\lambda }\right) \in \overline{B_{0}\left( 2R\right) }$\emph{\ be
the unique minimizer on the set }$\overline{B_{0}\left( 2R\right) }$ \emph{\
of the functional }$F_{\lambda ,\alpha }\left( r,s\right) $\emph{, which was
found in by Theorem 5.1, }%
\begin{equation}
F_{\lambda ,\alpha }\left( \widetilde{r}_{\min ,\lambda ,\alpha },\widetilde{
s}_{\min ,\lambda ,\alpha }\right) =\min_{\overline{B_{0}\left( 2R\right) }
}F_{\lambda ,\alpha }\left( r,s\right) .  \label{5.130}
\end{equation}%
\emph{\ Denote }%
\begin{equation}
\left( \overline{r}_{\min ,\lambda ,\alpha },\overline{s}_{\min ,\lambda
,\alpha }\right) =\left( \widetilde{r}_{\min ,\lambda ,\alpha }+P_{1}, 
\widetilde{s}_{\min ,\lambda ,\alpha }+P_{2}\right) .  \label{5.14}
\end{equation}%
\emph{Let the regularization parameter }%
\begin{equation}
\alpha =\xi ^{2}.  \label{5.140}
\end{equation}%
\emph{Then}%
\[
\left( \overline{r}_{\min ,\lambda ,\alpha },\overline{s}_{\min ,\lambda
,\alpha }\right) \in \overline{B\left( 3R\right) } 
\]%
\emph{and the following accuracy estimates hold: }%
\begin{equation}
\left. 
\begin{array}{c}
\left\Vert \left( \overline{r}_{\min ,\lambda ,\alpha },\overline{s}_{\min
,\lambda ,\alpha }\right) -\left( r^{\ast },s^{\ast }\right) \right\Vert
_{H^{2}\left( \Omega \right) \times H^{2}\left( \Omega \right) }\leq
C_{1}\xi \exp \left( \lambda \left( a+c_{1}\right) ^{2}\right) , \\ 
\forall \varphi \in \left( 0,2\pi \right) .%
\end{array}
\right.  \label{5.141}
\end{equation}%
\emph{\ }%
\begin{equation}
\left\Vert a_{\mbox{min},\lambda ,\alpha }-a^{\ast }\right\Vert
_{L_{2}\left( \Omega \right) }\leq C_{1}\xi \exp \left( \lambda \left(
a+c_{1}\right) ^{2}\right) ,  \label{5.15}
\end{equation}%
\emph{where the function }$a_{\mbox{min},\lambda }\left( \mathbf{x}\right) $%
\emph{\ is found via (\ref{3.25}) and (\ref{4.5}). }

\subsection{Global convergence of the gradient descent method}

\label{sec:5.3}

Let $\lambda _{2}$ be the number defined in (\ref{5.13}). Since $\lambda
_{2}=\lambda _{2}\left( R,\Omega ,\varepsilon \right) $ depends on the same
parameters as the number $C_{1},$ then, using (\ref{5.141}) and (\ref{5.15}%
), we denote%
\begin{equation}
\overline{C}_{1}=\overline{C}_{1}\left( R,\Omega ,\varepsilon \right)
=C_{1}\exp \left( \lambda _{2}\left( a+c_{1}\right) ^{2}\right)  \label{1}
\end{equation}%
and use notation (\ref{1}) in this subsection for different positive
constants depending only on these parameters.

\bigskip Similarly with (\ref{5.8}) we assume now that 
\begin{equation}
\left. 
\begin{array}{c}
R/3-\overline{C}_{1}\xi >0, \\ 
\left( r^{\ast },s^{\ast }\right) \left( \mathbf{x},\varphi \right) \in
B^{\ast }\left( R/3-\overline{C}_{1}\xi \right) .%
\end{array}
\right.  \label{5.16}
\end{equation}%
Assume that (\ref{5.140}) holds. Recalling $\left( \overline{r}_{\min
,\lambda _{2},\alpha },\overline{s}_{\min ,\lambda _{2},\alpha }\right) $ is
the pair of functions defined in (\ref{5.14}). It follows from (\ref{4.1}), (%
\ref{4.2}), (\ref{5.7}), (\ref{5.141}), (\ref{1}) and (\ref{5.16}) that it
is reasonable to assume that%
\begin{equation}
\left\Vert \left( \overline{r}_{\min ,\lambda ,\alpha },\overline{s}_{\min
,\lambda ,\alpha }\right) -\left( r^{\ast },s^{\ast }\right) \right\Vert
_{H_{2}}\leq \overline{C}_{1}\xi .  \label{2}
\end{equation}%
Using (\ref{5.16}) and (\ref{2}), we obtain 
\begin{equation}
\left( \overline{r}_{\min ,\lambda _{2},\alpha },\overline{s}_{\min ,\lambda
_{2},\alpha }\right) \in B\left( \frac{R}{3}\right) .  \label{5.160}
\end{equation}

Let the number $\gamma \in \left( 0,1\right) $ and let%
\begin{equation}
\left( r_{0},s_{0}\right) \left( \mathbf{x},\varphi \right) \in B\left( 
\frac{R}{3}\right)  \label{5.17}
\end{equation}%
be an arbitrary point of the set $B\left( R/3\right) $. We construct the
gradient descent method as:%
\begin{equation}
\hspace{-1.5 cm}
\left( r_{n},s_{n}\right) \left( \mathbf{x},\varphi \right) =\left(
r_{n-1},s_{n-1}\right) \left( \mathbf{x},\varphi \right) -\gamma J_{\lambda
_{2},\alpha }^{\prime }\left( r_{n-1},s_{n-1}\right) \left( \mathbf{x}
,\varphi \right) ,\mbox{ }n=1,2,...  \label{5.18}
\end{equation}%
Note that since by Theorem 5.1 $J_{\lambda _{2},\alpha }^{\prime }\left(
r_{n-1},s_{n-1}\right) $ satisfies (\ref{5.01}), then all terms of sequence (%
\ref{5.18}) have the same boundary conditions (\ref{3.21})-(\ref{3.24}).

\textbf{Theorem 5.3}. \emph{Let conditions (\ref{5.13}), (\ref{5.140}), (\ref%
{5.16})-(\ref{5.18}) hold}$.$ \emph{Then there exists a sufficiently small
number }$\gamma \in \left( 0,1\right) $\emph{\ and a number }$\theta =\theta
\left( \gamma \right) \in \left( 0,1\right) $\emph{\ such that all terms of
sequence (\ref{5.18}) belong to }$B\left( R\right) $\emph{\ and the
following estimates hold:}%
\begin{equation}
\left. 
\begin{array}{c}
\left\Vert \left( r_{n},s_{n}\right) -\left( r^{\ast },s^{\ast }\right)
\right\Vert _{H_{2}}\leq \\ 
\leq \overline{C}_{1}\xi +\theta ^{n}\left\Vert \left( r_{\min ,\lambda
_{2},\alpha },s_{\min ,\lambda _{2},\alpha }\right) -\left(
r_{0},s_{0}\right) \right\Vert _{H_{2}},%
\end{array}
\right.  \label{5.19}
\end{equation}%
\begin{equation}
\left. 
\begin{array}{c}
\left\Vert a_{n,\mbox{min},\lambda _{2},\alpha }-a^{\ast }\right\Vert
_{L_{2}\left( \Omega \right) }\leq \\ 
\leq \overline{C}_{1}\xi +\theta ^{n}\left\Vert \left( r_{\min ,\lambda
_{2},\alpha },s_{\min ,\lambda _{2},\alpha }\right) -\left(
r_{0},s_{0}\right) \right\Vert _{H_{2}},%
\end{array}
\right.  \label{5.20}
\end{equation}%
\emph{where functions }$a_{n,\mbox{min},\lambda _{2},\alpha }$\emph{\ are
found via (\ref{3.25}) and (\ref{4.5}) with the replacement in (\ref{4.5})
of }$\left( r_{\min },s_{\min }\right) $\emph{\ with }$\left(
r_{n},s_{n}\right) .$

\textbf{Proof.} Assuming that Theorems 5.1 and 5.2 are valid the proof
follows immediately from Theorem 6 of \cite{SAR}. $\square $

\textbf{Remark 5.1.}\emph{\ Definition of the global convergence given in
section 1 as well as (\ref{5.17}) imply that actually Theorem 5.3 claims the
global convergence. This is because the number }$R$\emph{\ is not assumed to
be small. }

\section{Proof of Theorem 5.1}

\label{sec:6}

Let $\left( r_{1},s_{1}\right) ,\left( r_{2},s_{2}\right) \in \overline{%
B\left( R\right) }$ be two arbitrary pairs of functions. Consider their
difference 
\begin{equation}
\left( r_{2},s_{2}\right) -\left( r_{1},s_{1}\right) =\left(
h_{1},h_{2}\right) \in \overline{B_{0}\left( 2R\right) }.  \label{6.1}
\end{equation}%
Then the pair of functions $\left( h_{1},h_{2}\right) \left( \mathbf{x}%
,\varphi \right) $ holds the same properties as the ones in (\ref{4.02}) for
the pair $\left( r,s\right) \left( \mathbf{x},\varphi \right) .$ By (\ref%
{3.19}) and (\ref{6.1}) we have for $\mathbf{x}\in \Omega ,\varphi \in
\left( 0,2\pi \right) $%
\[
L_{1}\left( r_{2},s_{2}\right) =L_{1}\left( r_{1}+h_{1},s_{1}+h_{2}\right)
=\Delta r_{1}+\Delta h_{1}+ 
\]%
\[
+2\nabla \left( r_{1}+h_{1}\right) \nabla \left( \frac{\left(
r_{1}-s_{1}\right) +\left( h_{1}-h_{2}\right) }{\varepsilon }\right) . 
\]%
We now single out the linear, with respect to $\left( h_{1},h_{2}\right) ,$
part of this expression. We have 
\[
L_{1}\left( r_{2},s_{2}\right) =L_{1}\left( r_{1},s_{1}\right) + 
\]%
\begin{equation}
+\Delta h_{1}+2\nabla r_{1}\nabla \left( \frac{h_{1}-h_{2}}{\varepsilon }%
\right) +2\nabla h_{1}\nabla \left( \frac{r_{1}-s_{1}}{\varepsilon }\right) +
\label{6.2}
\end{equation}%
\[
+2\nabla h_{1}\nabla \left( \frac{h_{1}-h_{2}}{\varepsilon }\right) . 
\]%
Hence, the linear with respect to $\left( h_{1},h_{2}\right) ,$ part of the
right hand side of (\ref{6.2}) is 
\begin{equation}
\hspace{-1 cm}
L_{1,\mbox{lin}}\left( h_{1},h_{2}\right) =\Delta h_{1}+2\nabla r_{1}\nabla
\left( \frac{h_{1}-h_{2}}{\varepsilon }\right) +2\nabla h_{1}\nabla \left( 
\frac{r_{1}-s_{1}}{\varepsilon }\right) .  \label{6.3}
\end{equation}%
Using (\ref{6.2}) and (\ref{6.3}), we obtain 
\[
\left( L_{1}\left( r_{2},s_{2}\right) \right) ^{2}-\left( L_{1}\left(
r_{1},s_{1}\right) \right) ^{2}=2L_{1}\left( r_{1},s_{1}\right) L_{1,%
\mbox{			lin}}\left( h_{1},h_{2}\right) + 
\]%
\begin{equation}
+\left( L_{1,\mbox{lin}}\left( h_{1},h_{2}\right) \right) ^{2}+4L_{1,%
\mbox{			lin}}\left( h_{1},h_{2}\right) \nabla h_{1}\nabla \left( \frac{%
h_{1}-h_{2}}{\varepsilon }\right) +  \label{6.4}
\end{equation}%
\[
+2L_{1}\left( r_{1},s_{1}\right) \nabla h_{1}\nabla \left( \frac{h_{1}-h_{2}%
}{\varepsilon }\right) +\left[ \nabla h_{1}\nabla \left( \frac{h_{1}-h_{2}}{%
\varepsilon }\right) \right] ^{2}. 
\]%
Using (\ref{4.02}), (\ref{6.3}), (\ref{6.4}) and Cauchy-Schwarz inequality,
we obtain 
\[
\left( L_{1}\left( r_{2},s_{2}\right) \right) ^{2}-\left( L_{1}\left(
r_{1},s_{1}\right) \right) ^{2}=2L_{1}\left( r_{1},s_{1}\right) L_{1,%
\mbox{			lin}}\left( h_{1},h_{2}\right) \geq 
\]%
\begin{equation}
\geq C_{1}\left( \Delta h_{1}\right) ^{2}-C_{1}\left( \left\vert \nabla
h_{1}\right\vert ^{2}+\left\vert \nabla h_{2}\right\vert ^{2}\right) ,%
\mbox{ 
	}\forall \left( \mathbf{x},\varphi \right) \in \Omega \times \left( 0,2\pi
\right) .  \label{6.5}
\end{equation}%
A similar procedure applied to the operator $L_{2}$ in (\ref{3.20}) leads to
the inequality, which is similar with (\ref{6.5}), 
\[
\left( L_{2}\left( r_{2},s_{2}\right) \right) ^{2}-\left( L_{2}\left(
r_{1},s_{1}\right) \right) ^{2}=2L_{2}\left( r_{1},s_{1}\right) L_{2,%
\mbox{			lin}}\left( h_{1},h_{2}\right) \geq 
\]%
\begin{equation}
\geq C_{1}\left( \Delta h_{2}\right) ^{2}-C_{1}\left( \left\vert \nabla
h_{1}\right\vert ^{2}+\left\vert \nabla h_{2}\right\vert ^{2}\right) ,%
\mbox{ 
	}\forall \left( \mathbf{x},\varphi \right) \in \Omega \times \left( 0,2\pi
\right) ,  \label{6.6}
\end{equation}%
where $L_{2,\mbox{lin}}\left( h_{1},h_{2}\right) $ is linear with respect to 
$\left( h_{1},h_{2}\right) .$

Denote 
\[
\left[ u,v\right] ,\forall u,v\in H^{3}\left( \Omega \right) \times
H^{3}\left( \Omega \right) . 
\]%
the scalar product in the Hilbert space $H^{3}\left( \Omega \right) \times
H^{3}\left( \Omega \right) .$ Using (\ref{4.4}), (\ref{6.5}) and the obvious
analog of (\ref{6.5}) for the operator $L_{2},$ we obtain for all $\varphi
\in \left( 0,2\pi \right) $ 
\[
J_{\lambda ,\alpha }\left( r_{2},s_{2}\right) \left( \varphi \right)
-J_{\lambda ,\alpha }\left( r_{1},s_{1}\right) \left( \varphi \right) = 
\]%
\[
\hspace{-1.5cm}=2\sqrt{\varepsilon }\int\limits_{\Omega }\left[ L_{1}\left(
r_{1},s_{1}\right) L_{1,\mbox{lin}}\left( h_{1},h_{2}\right) +L_{2}\left(
r_{1},s_{1}\right) L_{2,\mbox{lin}}\left( h_{1},h_{2}\right) \right] \psi
_{\lambda }\left( \mathbf{x}\right) d\mathbf{x}+ 
\]%
\[
+2\alpha \left[ \left( r_{1},s_{1}\right) ,\left( h_{1},h_{2}\right) \right]
+ 
\]%
\begin{equation}
\hspace{-1.5cm}+\sqrt{\varepsilon }\int\limits_{\Omega }\left[ \left( L_{1,%
\mbox{lin}}\left( h_{1},h_{2}\right) \right) ^{2}+4L_{1,\mbox{lin}}\left(
h_{1},h_{2}\right) \nabla h_{1}\nabla \left( \frac{h_{1}-h_{2}}{\varepsilon }%
\right) \right] \psi _{\lambda }\left( \mathbf{x}\right) d\mathbf{x}+
\label{6.9}
\end{equation}%
\[
\hspace{-1.5cm}+\sqrt{\varepsilon }\int\limits_{\Omega }\left[ \left\{
2L_{1}\left( r_{1},s_{1}\right) \nabla h_{1}\nabla \left( \frac{h_{1}-h_{2}}{%
\varepsilon }\right) +\left[ \nabla h_{1}\nabla \left( \frac{h_{1}-h_{2}}{%
\varepsilon }\right) \right] ^{2}\right\} \right] \psi _{\lambda }\left( 
\mathbf{x}\right) d\mathbf{x+} 
\]%
\[
+S\left( h_{1},h_{2}\right) +\alpha \left\Vert \left( h_{1},h_{2}\right)
\right\Vert _{H^{3}\left( \Omega \right) \times H^{3}\left( \Omega \right)
}^{2},\mbox{ }\varphi \in \left( 0,2\pi \right) . 
\]%
where the term $S\left( h_{1},h_{2}\right) $ is similar with the ones in the
lines 4 and 5\textbf{\ of }(\ref{6.9}), except that it is generated by the
operator $L_{2}$ in (\ref{3.20}) rather than by the operator $L_{1}$ in (\ref%
{3.19}). Consider the expression in the second and third lines of (\ref{6.9}%
),%
\begin{equation}
\hspace{-1.5cm}\left. 
\begin{array}{c}
\widehat{J}_{\lambda ,\alpha ,r_{1},s_{1}}\left( h_{1},h_{2}\right) \left(
\varphi \right) = \\ 
=2\sqrt{\varepsilon }\int\limits_{\Omega }\left[ L_{1}\left(
r_{1},s_{1}\right) L_{1,\mbox{lin}}\left( h_{1},h_{2}\right) +L_{2}\left(
r_{1},s_{1}\right) L_{2,\mbox{lin}}\left( h_{1},h_{2}\right) \right] \psi
_{\lambda }\left( \mathbf{x}\right) d\mathbf{x}+ \\ 
++2\alpha \left[ \left( r_{1},s_{1}\right) ,\left( h_{1},h_{2}\right) \right]
,\mbox{ }\varphi \in \left( 0,2\pi \right) .%
\end{array}%
\right.  \label{6.10}
\end{equation}%
It is clear that 
\begin{equation}
\widehat{J}_{\lambda ,\alpha ,r_{1},s_{1}}\left( h_{1},h_{2}\right) \left(
\varphi \right) :H_{0}^{3}\left( \Omega \right) \times H_{0}^{3}\left(
\Omega \right) \rightarrow \mathbb{R},\mbox{ }\forall \varphi \in \left(
0,2\pi \right)  \label{6.11}
\end{equation}%
is a bounded linear functional of $\left( h_{1},h_{2}\right) $ for every $%
\varphi \in \left( 0,2\pi \right) $. Hence, Riesz theorem implies that for
each $\varphi \in \left( 0,2\pi \right) $ there exists unique vector
function $\widetilde{J}_{\lambda ,\alpha ,r_{1},s_{1}}\left( \varphi \right)
\in H_{0}^{3}\left( \Omega \right) \times H_{0}^{3}\left( \Omega \right) $
such that 
\begin{equation}
\left. 
\begin{array}{c}
\widehat{J}_{\lambda ,\alpha ,r_{1},s_{1}}\left( h_{1},h_{2}\right) \left(
\varphi \right) =\left[ \widetilde{J}_{\lambda ,\alpha ,r_{1},s_{1}},\left(
h_{1},h_{2}\right) \right] \left( \varphi \right) , \\ 
\forall \left( h_{1},h_{2}\right) \in H_{0}^{3}\left( \Omega \right) \times
H_{0}^{3}\left( \Omega \right) .%
\end{array}%
\right.  \label{6.12}
\end{equation}%
It follows from (\ref{6.1})-(\ref{6.4}) and (\ref{6.9})-(\ref{6.12}) that
for all $\varphi \in \left( 0,2\pi \right) $%
\begin{equation}
\hspace{-1.5cm}\frac{J_{\lambda ,\alpha }\left(
r_{1}+h_{1},s_{1}+h_{2}\right) \left( \varphi \right) -J_{\lambda ,\alpha
}\left( r_{1},s_{1}\right) \left( \varphi \right) -\widehat{J}_{\lambda
,\alpha ,r_{1},s_{1}}\left( h_{1},h_{2}\right) \left( \varphi \right) }{%
\left\Vert \left( h_{1},h_{2}\right) \right\Vert _{H^{3}\left( \Omega
\right) \times H^{3}\left( \Omega \right) }}\rightarrow 0  \label{6.13}
\end{equation}%
as $\left\Vert \left( h_{1},h_{2}\right) \right\Vert _{H^{3}\left( \Omega
\right) \times H^{3}\left( \Omega \right) }\rightarrow 0.$ Hence, using (\ref%
{6.10})-(\ref{6.13}), we obtain that $\widetilde{J}_{\lambda ,\alpha
,r_{1},s_{1}}$ is the Fr\'{e}chet derivative of the functional $J_{\lambda
,\alpha }$ at the point $\left( r_{1},s_{1}\right) ,$ i.e. 
\begin{equation}
\widetilde{J}_{\lambda ,q_{1},r_{1},s_{1}}\left( \varphi \right) =J_{\lambda
,\alpha }^{\prime }\left( r_{1},s_{1}\right) \left( \varphi \right) \in
H_{0}^{3}\left( \Omega \right) \times H_{0}^{3}\left( \Omega \right) ,%
\mbox{ 
	}\forall \varphi \in \left( 0,2\pi \right) .  \label{6.14}
\end{equation}%
We omit the proof of the Lipschitz continuity property (\ref{5.2}) of $%
J_{\lambda ,\alpha }^{\prime }\left( r,s\right) \left( \varphi \right) $
since this proof is similar with the proof of Theorem 5.3.1 of \cite{KL}.

We now prove the strong convexity property (\ref{5.3}). Using Theorem 5.1, (%
\ref{6.5}), (\ref{6.6}), (\ref{6.9}), (\ref{6.12}) and (\ref{6.14}), we
obtain%
\[
\hspace{-1.5 cm} J_{\lambda ,\alpha }\left( r_{1}+h_{1},s_{1}+h_{2}\right)
\left( \varphi \right) -J_{\lambda ,\alpha }\left( r_{1},s_{1}\right) \left(
\varphi \right) -\left[ J_{\lambda ,\alpha }^{\prime }\left(
r_{1},s_{1}\right) \left( \varphi \right) ,\left( h_{1},h_{2}\right) \right]
\geq 
\]%
\[
\hspace{-1.5 cm} \geq C_{1}\int\limits_{\Omega }\left[ \left( \Delta
h_{1}\right) ^{2}+\left( \Delta h_{2}\right) ^{2}\right] \psi _{\lambda
}\left( \mathbf{x} \right) d\mathbf{x}-C_{1}\int\limits_{\Omega }\left[
\left\vert \nabla h_{1}\right\vert ^{2}+\left\vert \nabla h_{2}\right\vert
^{2}\right] \psi _{\lambda }\left( \mathbf{x}\right) d\mathbf{x}\geq 
\]%
\begin{equation}
\geq \frac{C_{1}}{\lambda }\int\limits_{\Omega }\left(
h_{1xx}^{2}+h_{1xy}^{2}+h_{1yy}^{2}+h_{2xx}^{2}+h_{2xy}^{2}+h_{2yy}^{2}
\right) \psi _{\lambda }\left( \mathbf{x}\right) d\mathbf{x}+  \label{6.15}
\end{equation}%
\[
+C_{1}\lambda \int\limits_{\Omega }\left[ \left( \nabla h_{1}\right)
^{2}+\left( \nabla h_{2}\right) ^{2}+\lambda ^{2}\left(
h_{1}^{2}+h_{2}^{2}\right) \right] \psi _{\lambda }\left( \mathbf{x}\right)
d \mathbf{x}- 
\]%
\[
-C_{1}\int\limits_{\Omega }\left[ \left\vert \nabla h_{1}\right\vert
^{2}+\left\vert \nabla h_{2}\right\vert ^{2}\right] \psi _{\lambda }\left( 
\mathbf{x}\right) d\mathbf{x}+\alpha \left\Vert \left( h_{1},h_{2}\right)
\right\Vert _{H^{3}\left( \Omega \right) \times H^{3}\left( \Omega \right)
}^{2},\mbox{ } 
\]%
\[
\forall \lambda \geq \lambda _{0},\mbox{ }\forall \varphi \in \left( 0,2\pi
\right) . 
\]%
Hence, we can choose a sufficiently large number $\lambda _{1}=\lambda
_{1}\left( R,\Omega ,\varepsilon \right) \geq \lambda _{0}$ such that (\ref%
{6.15}) becomes%
\[
\hspace{-1.5 cm} J_{\lambda ,\alpha }\left( r_{1}+h_{1},s_{1}+h_{2}\right)
\left( \varphi \right) -J_{\lambda ,\alpha }\left( r_{1},s_{1}\right) \left(
\varphi \right) -\left[ J_{\lambda ,\alpha }^{\prime }\left(
r_{1},s_{1}\right) \left( \varphi \right) ,\left( h_{1},h_{2}\right) \right]
\geq 
\]%
\begin{equation}
\hspace{-1 cm} \geq C_{1}\exp \left( \lambda \left( a-c_{1}\right) \right)
\left\Vert \left( h_{1},h_{2}\right) \right\Vert _{H^{2}\left( \Omega
\right) \times H^{2}\left( \Omega \right) }^{2}+\alpha \left\Vert \left(
h_{1},h_{2}\right) \right\Vert _{H^{3}\left( \Omega \right) \times
H^{3}\left( \Omega \right) }^{2},  \label{6.16}
\end{equation}%
\[
\forall \lambda \geq \lambda _{1},\mbox{ }\forall \varphi \in \left( 0,2\pi
\right) . 
\]%
which proves (\ref{5.3}).

Existence and uniqueness of the minimizer $\left( r_{\min ,\lambda ,\alpha
},s_{\min ,\lambda ,\alpha }\right) \left( \mathbf{x},\varphi \right) \in 
\overline{B\left( R\right) }$\ of the functional $J_{\lambda ,\alpha }\left(
r,s\right) \left( \varphi \right) $\ on the set $\overline{B\left( R\right) }
$ as well as inequality (\ref{5.4}) follow immediately from (\ref{5.3}) and
a combination of Lemma 5.2.1 with Theorem 5.2.1 of \cite{KL}. $\square $

\section{Proof of Theorem 5.2}

\label{sec:7}

In this section $\lambda \geq \lambda _{2},$ where $\lambda _{2}$ is defined
in (\ref{5.13}). Let $F_{\lambda ,\alpha }\left( r,s\right) $ be the
functional defined in (\ref{5.12}). Recall that an obvious analog of Theorem
5.1 is valid for $F_{\lambda ,\alpha }\left( r,s\right) $. In addition,
recall that $\left( \widetilde{r}_{\min ,\lambda ,\alpha },\widetilde{s}%
_{\min ,\lambda ,\alpha }\right) \in \overline{B_{0}\left( 2R\right) }$\emph{%
\ \ }is the unique minimizer on the set $\overline{B_{0}\left( 2R\right) }$
of the functional $F_{\lambda ,\alpha }\left( r,s\right) $ for $\lambda .$
By (\ref{5.3}), (\ref{5.10}) and the second line of (\ref{5.12}) we have for 
$\varphi \in \left( 0,2\pi \right) :$%
\[
F_{\lambda ,\alpha }\left( \widetilde{r}^{\ast },\widetilde{s}^{\ast
}\right) \left( \varphi \right) -F_{\lambda ,\alpha }\left( \widetilde{r}
_{\min ,\lambda ,\alpha },\widetilde{s}_{\min ,\lambda ,\alpha }\right)
\left( \varphi \right) - 
\]%
\begin{equation}
-\left[ F_{\lambda ,\alpha }^{\prime }\left( \widetilde{r}_{\min ,\lambda
,\alpha },\widetilde{s}_{\min ,\lambda ,\alpha }\right) \left( \varphi
\right) ,\left( \widetilde{r}^{\ast }-\widetilde{r}_{\min ,\lambda ,\alpha
}, \widetilde{s}^{\ast }-\widetilde{s}_{\min ,\lambda ,\alpha }\right)
\left( \varphi \right) \right] \geq  \label{7.1}
\end{equation}%
\[
\geq C_{1}\exp \left( \lambda \left( a-c_{1}\right) \right) \left\Vert
\left( \widetilde{r}^{\ast }-\widetilde{r}_{\min ,\lambda ,\alpha }, 
\widetilde{s}^{\ast }-\widetilde{s}_{\min ,\lambda ,\alpha }\right) \left(
\varphi \right) \right\Vert _{H^{2}\left( \Omega \right) \times H^{2}\left(
\Omega \right) }^{2}. 
\]%
By (\ref{5.4})%
\[
-\left[ F_{\lambda ,\alpha }^{\prime }\left( \widetilde{r}_{\min ,\lambda
,\alpha },\widetilde{s}_{\min ,\lambda ,\alpha }\right) \left( \varphi
\right) ,\left( \widetilde{r}^{\ast }-\widetilde{r}_{\min ,\lambda ,\alpha
}, \widetilde{s}^{\ast }-\widetilde{s}_{\min ,\lambda ,\alpha }\right)
\left( \varphi \right) \right] \leq 0. 
\]%
Also, obviously $-F_{\lambda ,\alpha }\left( \widetilde{r}_{\min ,\lambda
,\alpha },\widetilde{s}_{\min ,\lambda ,\alpha }\right) \left( \varphi
\right) \leq 0.$ Hence, (\ref{7.1}) implies for $\varphi \in \left( 0,2\pi
\right) $ 
\[
F_{\lambda ,\alpha }\left( \widetilde{r}^{\ast },\widetilde{s}^{\ast
}\right) \left( \varphi \right) \geq 
\]%
\begin{equation}
\geq C_{1}\exp \left( \lambda \left( a-c_{1}\right) \right) \left\Vert
\left( \widetilde{r}^{\ast }-\widetilde{r}_{\min ,\lambda ,\alpha }, 
\widetilde{s}^{\ast }-\widetilde{s}_{\min ,\lambda ,\alpha }\right) \left(
\varphi \right) \right\Vert _{H^{2}\left( \Omega \right) \times H^{2}\left(
\Omega \right) }^{2}.  \label{7.2}
\end{equation}%
Next, by (\ref{5.10}) and (\ref{5.12}) 
\[
F_{\lambda ,\alpha }\left( \widetilde{r}^{\ast },\widetilde{s}^{\ast
}\right) \left( \varphi \right) =J_{\lambda ,\alpha }\left( \widetilde{r}
^{\ast }+P_{1},\widetilde{s}^{\ast }+P_{2}\right) \left( \varphi \right) = 
\]%
\begin{equation}
=J_{\lambda ,\alpha }\left( \left( \widetilde{r}^{\ast }+P_{1}^{\ast
}\right) +\left( P_{1}-P_{1}^{\ast }\right) ,\left( \widetilde{s}^{\ast
}+P_{2}^{\ast }\right) +\left( P_{2}-P_{2}^{\ast }\right) \right) \left(
\varphi \right) =  \label{7.3}
\end{equation}%
\[
=J_{\lambda ,\alpha }\left( r^{\ast }+\left( P_{1}-P_{1}^{\ast }\right)
,s^{\ast }+\left( P_{2}-P_{2}^{\ast }\right) \right) . 
\]%
Next, by (\ref{4.4}) and (\ref{5.6}) and (\ref{5.140}) and 
\begin{equation}
\hspace{-1.5 cm} \left. 
\begin{array}{c}
J_{\lambda ,\alpha }\left( r^{\ast }+\left( P_{1}-P_{1}^{\ast }\right)
,s^{\ast }+\left( P_{2}-P_{2}^{\ast }\right) \right) \left( \varphi \right) =
\\ 
=\sqrt{\varepsilon }\int\limits_{\Omega }\left( L_{1}\left( r^{\ast }+\left(
P_{1}-P_{1}^{\ast }\right) ,s^{\ast }+\left( P_{2}-P_{2}^{\ast }\right)
\right) \left( \mathbf{x},\varphi \right) \right) ^{2}\psi _{\lambda }\left( 
\mathbf{x}\right) d\mathbf{x}+ \\ 
+\sqrt{\varepsilon }\int\limits_{\Omega }\left( L_{2}\left( r^{\ast }+\left(
P_{1}-P_{1}^{\ast }\right) ,s^{\ast }+\left( P_{2}-P_{2}^{\ast }\right)
\right) \left( \mathbf{x},\varphi \right) \right) ^{2}\psi _{\lambda }\left( 
\mathbf{x}\right) d\mathbf{x+} \\ 
+\xi ^{2}\left\Vert r^{\ast }+\left( P_{1}-P_{1}^{\ast }\right) ,s^{\ast
}+\left( P_{2}-P_{2}^{\ast }\right) \left( \mathbf{x},\varphi \right)
\right\Vert _{H^{3}\left( \Omega \right) \times H^{3}\left( \Omega \right)
}^{2},\mbox{ }\forall \varphi \in \left( 0,2\pi \right) .%
\end{array}
\right.  \label{7.4}
\end{equation}

Obviously $L_{1}\left( r^{\ast },s^{\ast }\right) =L_{2}\left( r^{\ast
},s^{\ast }\right) =0.$ Hence, (\ref{2.20}), (\ref{4.3}), (\ref{5.5}), (\ref%
{5.140}), (\ref{7.3}) and (\ref{7.4}) imply%
\[
F_{\lambda ,\alpha }\left( \widetilde{r}^{\ast },\widetilde{s}^{\ast
}\right) \left( \varphi \right) \leq C_{1}\exp \left( 2\lambda \left(
a+c_{1}\right) \right) \xi ^{2}. 
\]%
Combining this with (\ref{5.14}) and (\ref{7.2}), we obtain (\ref{5.141}).
Estimate (\ref{5.15}) follows from (\ref{3.25}), (\ref{4.5}) and (\ref{5.141}%
). \ $\square $

\section{Numerical Studies}

\bigskip \label{sec:8}

\begin{figure}[tbph]
\centering
\includegraphics[width = 3in]{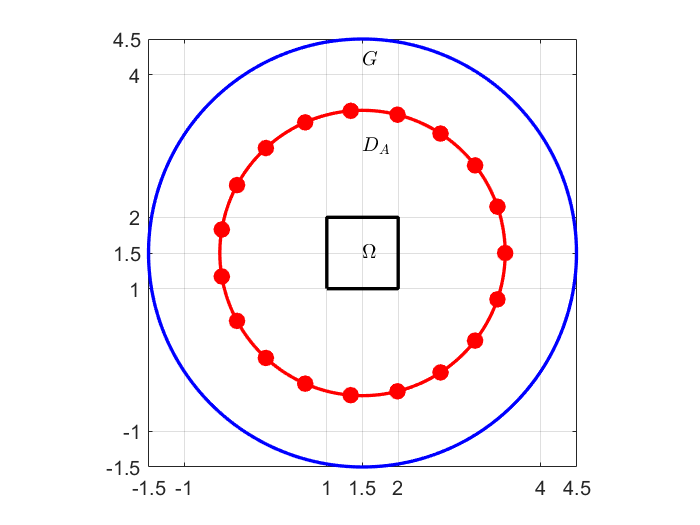}
\caption{A schematic diagram of our measurement setup. The small black
square indicates the domain $\Omega $ defined in (\protect\ref{2.20}). The
large blue circle is the domain $G$ where the forward problem (\protect\ref%
{2.8}), (\protect\ref{2.9}) is solved to generate the data for the inverse
problem. The red circle is the circle $C_{A}$ in (\protect\ref{2.1}) and
small red discs on it indicate positions of the point source.}
\label{plot_domain}
\end{figure}

For the disk $D_{A}\left( a,b\right) $ in (\ref{2.1}), we choose $%
a=b=1.5,A=2 $. For the domain $\Omega $ in (\ref{2.20}), we choose $%
c_{1}=c_{2}=0.5$. In the forward problem (\ref{2.8}), (\ref{2.9}), we choose 
$G$ as the disk of the radius 3 with the center at $(1.5,1.5)$. We display
the corresponding schematic diagram in Figure \ref{plot_domain}.

We choose $\rho =0.1$ in (\ref{2.3}) for the function $f\left( \mathbf{x}-%
\mathbf{x}_{0}\right) $. For the conductivity coefficient $\sigma (\mathbf{x}%
)$ to be reconstructed, we set 
\begin{equation}
\sigma (\mathbf{x})=\left\{ 
\begin{array}{cc}
\sigma_{a}=const.>1, & \mbox{inside the tested inclusion,} \\ 
1, & \mbox{outside the tested inclusion.}%
\end{array}%
\right.  \label{8.01}
\end{equation}%
To make sure $\sigma \left( \mathbf{x}\right) \in C^{1}\left( \overline{G}%
\right) $, we slightly smooth out $\sigma \left( \mathbf{x}\right) $ near
the boundaries of our tested inclusions. Then we set: 
\begin{eqnarray}
&&\hspace{-0.5cm}\mbox{correct inclusion/background contrast}=\frac{\sigma
_{a}}{1},  \label{8.02} \\
&&\hspace{-1cm}\mbox{computed inclusion/background contrast}=\frac{\max_{%
\mbox{inclusion}}\left( \sigma _{\mbox{comp}}(\mathbf{x})\right) }{1},
\label{8.03}
\end{eqnarray}%
where $\sigma _{\mbox{comp}}(\mathbf{x})$ is the computed coefficient $%
\sigma (\mathbf{x}).$ In the numerical tests below, we take $\sigma _{a}=2,4$%
,$8$ which correspond to 2:1, 4:1 and 8:1 inclusion/background contrasts
respectively. In the first series of numerical experiments we test the
inclusions with the shapes of the letters `$A$' and `$\Omega $'. In the
second series we test the shapes of CT scans of an abdomen.

\begin{figure}[tbph]
\centering
\includegraphics[width = 7in]{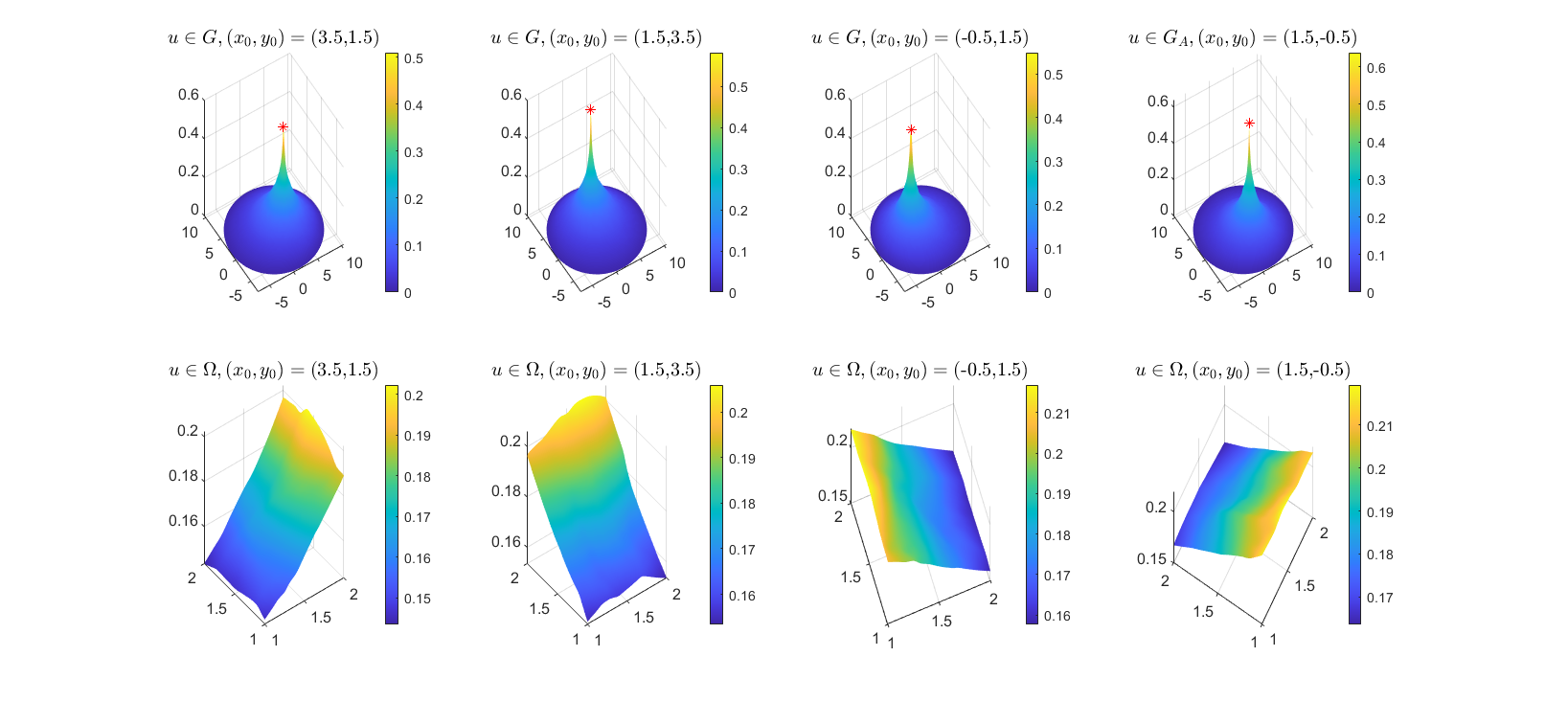}
\caption{The results of the solution of the forward problem (\protect\ref%
{2.8}), (\protect\ref{2.9}) for the case when the inclusion has the shape of
the letter `$A$' with $\protect\sigma _{a}=2$ inside of this letter, see (%
\protect\ref{8.01}). Top $\mathbf{x}\in G$, bottom $\mathbf{x}\in \Omega .$
The positions of the source are: $\mathbf{x}_{0} = (3.5, 1.5)$ (1st column), 
$\mathbf{x}_{0} = (1.5, 3.5)$ (2nd column), $\mathbf{x}_{0} = (-0.5, 1.5)$
(3rt column), $\mathbf{x}_{0} = (1.5, -0.5)$ (4th column). The red star is
the source position.}
\label{plot_A9}
\end{figure}

We choose the finite element method with $h=1/160$ to solve the forward
problem in (\ref{2.8})-(\ref{2.9}). For the source position $\mathbf{x}%
_{0}\in C_{A}$ in (\ref{3.11}), we choose $\varphi =nh_{\varphi }$, $%
n=1,2,\cdots ,199$, $h_{\varphi }=\pi /100$ . When the inclusion in (\ref%
{8.01}) has the shape of the letter `$A$' with $\sigma _{a}=2$, the results
of the solution of the forward problem of $\mathbf{x}\in G$ and for four
positions of the source 
\[
\mathbf{x}_{0}=(x_{0},y_{0})=(3.5,1.5),(1.5,3.5),(-0.5,1.5),(1.5,-0.5) 
\]%
are displayed in Figure \ref{plot_A9}.

To solve Coefficient Inverse Problem (\ref{2.11}), we have chosen the
spatial mesh sizes in the computations of the Minimization Problem 1 in (\ref%
{4.4}) and the Minimization Problem 2 in (\ref{4.9}) as $1/40\times 1/40$.
In both functionals $J_{\lambda ,\alpha }\left( r,s\right) \left( \varphi
\right) $ and $K\left( V\right) $ we write the differential operators in the
forms of finite differences and minimize the resulting discretized
functionals with respect to the values of corresponding functions at grid
points. Since $H^{3}\left( \Omega \right) -$norms are inconvenient for the
numerical implementation, we replace in (\ref{4.4}) the regularization term $%
\alpha \left\Vert \left( r,s\right) \left( \mathbf{x},\varphi \right)
\right\Vert _{H^{3}\left( \Omega \right) \times H^{3}\left( \Omega \right)
}^{2}$ with $\alpha \left\Vert \left( r,s\right) \left( \mathbf{x},\varphi
\right) \right\Vert _{H^{2}\left( \Omega \right) \times H^{2}\left( \Omega
\right) }^{2}.$ Since we actually work with not too many grid points and
since norms in finite dimensional spaces are equivalent, then this
replacement did have a negative impact on our numerical results. An
extension of the above theory on the discrete case is outside of the scope
of this publication.

To guarantee that the solution of the problem of the minimization of the
functional $J_{\lambda ,\alpha }\left( r,s\right) \left( \varphi \right) $
in (\ref{4.4}) satisfies the boundary conditions (\ref{3.21})-(\ref{3.24}),
we adopt the Matlab's built-in optimization toolbox \textbf{fmincon} to
minimize the discretized form of these corresponding functions. The same is
true for the Minimization Problem 2.

The iterations of \textbf{fmincon} stop when the condition 
\begin{equation}
|\nabla J_{\lambda ,\alpha }\left( r,s\right) \left( \varphi \right)
|<10^{-2},\quad |\nabla K(V)|<10^{-2}  \label{8.04}
\end{equation}%
are met. An explanation of the stopping criterion (\ref{8.04}) is provided
below.

We introduce the random noise in the observation data in \eref{2.11} as
follows: 
\begin{equation}
\left. 
\begin{array}{c}
g_{0}^{\xi _{0}}\left( \mathbf{x},\mathbf{x}_{0}\right) =g_{0}\left( \mathbf{%
x},\mathbf{x}_{0}\right) \left( 1+\delta \xi _{0}\left( \mathbf{x}\right)
\right) ,\quad \mathbf{x}\in \partial \Omega , \\ 
g_{1}^{\xi _{1}}\left( \mathbf{x},\mathbf{x}_{0}\right) =g_{1}\left( \mathbf{%
x},\mathbf{x}_{0}\right) \left( 1+\delta \xi _{1}\left( \mathbf{x}\right)
\right) ,\quad \mathbf{x}\in \Gamma ,%
\end{array}%
\right.  \label{8.05}
\end{equation}%
where $\xi _{0}$ is the uniformly distributed random variable in the
interval $[-1,1]$ depending on the point $\mathbf{x}\in \partial \Omega $.
Also, $\xi _{1}$ is the uniformly distributed random variable in the
interval $[-1,1]$ depending on the point $\mathbf{x}\in \Gamma $, and $%
\delta =0.03$ corresponds to the $3\%$ noise level. Since we deal with the
first $\varphi -$derivatives of the noisy functions $g_{0}^{\xi _{0}}\left( 
\mathbf{x},\mathbf{x}_{0}\right) $ and $g_{1}^{\xi _{1}}\left( \mathbf{x},%
\mathbf{x}_{0}\right) $, we have to design a numerical method to
differentiate the noisy data. First, we use the natural cubic splines to
approximate the noisy input data (\ref{8.05}). Next, we use the derivatives
of those splines to approximate the derivatives of corresponding noisy
observation data. We generate the corresponding cubic splines in $(0,2\pi )$
with the mesh grid size $h_{\varphi }=\pi /100$, and then we calculate their
derivatives to approximate the first derivatives with respect to $\varphi $.

We choose the optimal pair of parameters $\left( \alpha ,\varepsilon \right) 
$ by the trial and error procedure for the reference Test 1. For each
considered pair $\left( \alpha ,\varepsilon \right) ,$ we test different
values of the parameter $\lambda $ to obtain its optimal value $\lambda _{%
\mbox{opt}}\left( \alpha ,\varepsilon \right) $ for this pair. Once the so
chosen triple $\left( \alpha ,\varepsilon ,\lambda \right) $ of parameters
is selected, we consider it as the optimal choice of parameters. An
important point to make here is that exactly the same triple of optimal
parameters is used for all follow up tests when imaging letters below.
However, when using the CT scan of the abdomen below, we deal with a
different medium. This means that we repeat the procedure of our choice of
parameters again for this case.

\textbf{Remarks 8.1}:

\begin{enumerate}
\item \emph{As the test media, we intentionally choose letter-like shapes of
inclusions in the first series of numerical experiments and the CT scans of
the abdomen in the second series. This is done to demonstrate that our
technique works well for complicated media.}

\item \emph{The above procedure of the choice of an optimal triple }$\left(
\alpha ,\varepsilon ,\lambda \right) $\emph{\ of parameters is similar to
the conventional calibration procedure, which is often used in many real
World applications. Furthermore, quite similar procedures were used in all
above cited works \cite{SAR,KEIT,KL,Ktr,Ktransp,MFG} on the numerical
studies of the convexification method for CIPs.}

\item \emph{Even though theorems of our convergence analysis are valid only
for sufficiently large values of the parameter }$\lambda ,$\emph{\ we have
discovered in all our works on the convexification listed in item 2 that
actually optimal values of }$\lambda $\emph{\ belong to the interval }$%
\lambda \in \left[ 1,5\right] .$\emph{\ In fact, this is similar with many
asymptotic theories. Indeed, it is typically established in such a theory
that if a certain parameter }$X$\emph{\ is sufficiently large, then a
certain formula }$Y$\emph{\ provides a good approximation for a process.
However, it is also typical that in a computational practice only numerical
experiments can tell one which exactly values of }$X$\emph{\ are appropriate
ones. }
\end{enumerate}

\textbf{Test 1.} We test the case when the inclusion in (\ref{8.01}) has the
shape of the letter `$A$' with $\sigma _{a}=2$ inside of it. We use this
test as a reference one to figure out the optimal triple $\left( \alpha
,\varepsilon ,\lambda \right) $ of parameters. We have found that this
triple is 
\begin{equation}
\alpha =0.01,\varepsilon =0.0002,\lambda =3.  \label{8.6}
\end{equation}%
This is our optimal choice of parameters for the case when inclusions have
letter-like shapes.

\begin{figure}[tbph]
\centering
\includegraphics[width
=6in]{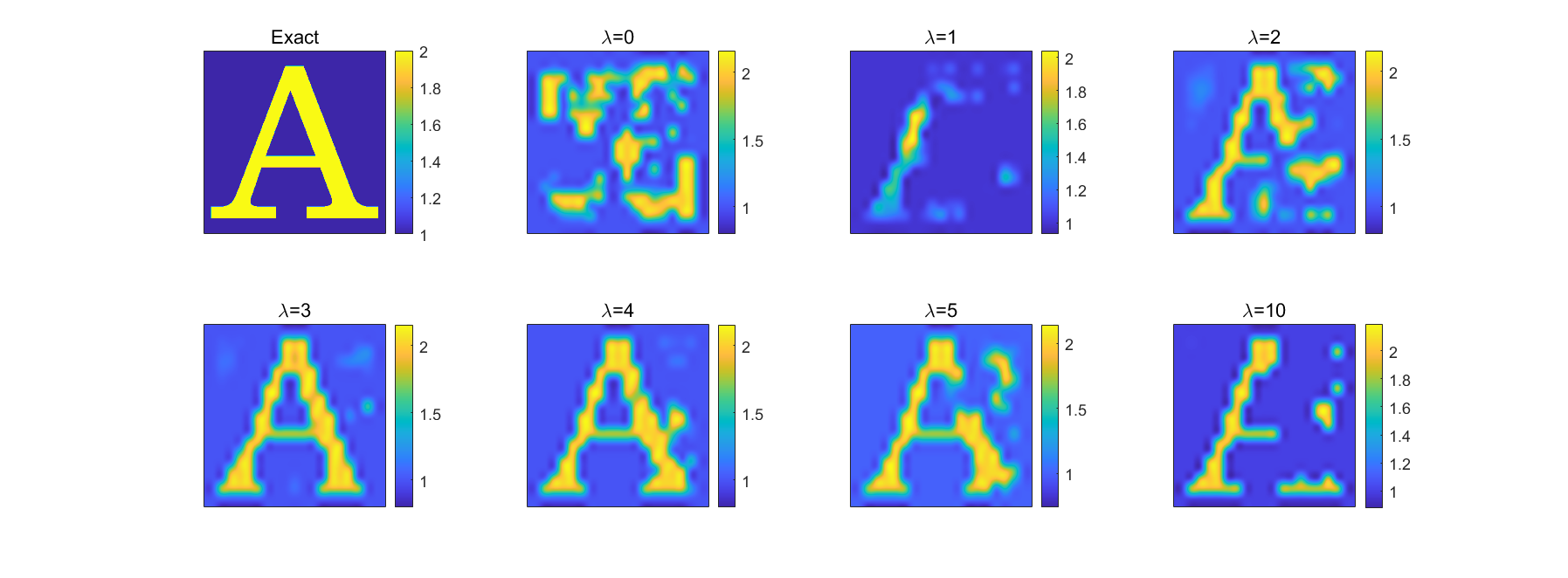}
\caption{Test 1. The reconstructed coefficient $\protect\sigma( \mathbf{x} ) 
$, where the function $\protect\sigma( \mathbf{x} ) $ is given in \eref{8.01}
with $\protect\sigma_{a}=2$ inside of the letter `A'. We test different
values of $\protect\lambda $.}
\label{plot_diff_lambda}
\end{figure}

Figure \ref{plot_diff_lambda} shows how do we choose the optimal value of $%
\lambda $ once the optimal pair $\left( \alpha ,\varepsilon \right) $ is
selected as in (\ref{8.6}). We observe that the images have a low quality
for $\lambda =0,1,2$. Then the quality is improved with $\lambda =3,4,5$,
and the reconstruction quality deteriorates at $\lambda =10$. On the other
hand, the image is accurate at $\lambda =3$, including the accurate
reconstruction of the inclusion/background contrast (\ref{8.02}).

\begin{figure}[tbph]
\centering
\includegraphics[width = 5in]{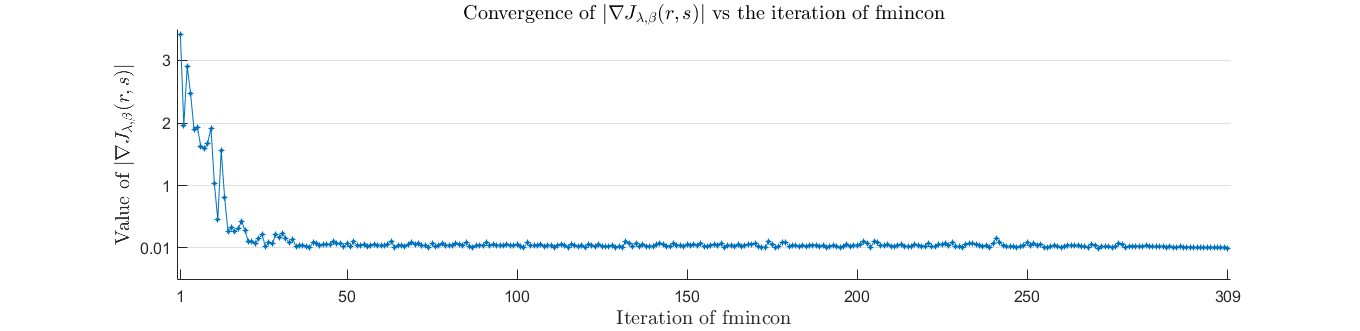}
\caption{Test 1. The convergence behavior of $\left\vert \protect\nabla J_{%
\protect\lambda ,\protect\alpha }\left( r,s\right) \left( \protect\varphi %
\right) \right\vert $ with the iterations of \textbf{fmincon} for $\protect%
\varphi =\protect\pi $ and the optimal triple of parameters $\left( \protect%
\alpha ,\protect\varepsilon ,\protect\lambda \right) $ as in (\protect\ref%
{8.6}). The function $\protect\sigma \left( \mathbf{x}\right) $ is given in (%
\protect\ref{8.01}) with $\protect\sigma _{a}=2$ inside of the letter `$A$'.
The $y-$axis corresponds to the value of $\left\vert \protect\nabla J_{%
\protect\lambda ,\protect\alpha }\left( r,s\right) \left( \protect\varphi %
\right) \right\vert$.}
\label{plot_grad}
\end{figure}

We display now in Figure \ref{plot_grad} the convergence behavior of $%
|\nabla J_{\lambda ,\alpha }(r,s)(\varphi )|$ with respect to the iterations
of \textbf{fmincon} for $\varphi =\pi $ and for the optimal triple $\left(
\alpha ,\varepsilon ,\lambda \right) $ of parameters as in (\ref{8.6}). This
figure explains the stopping criterion (\ref{8.04}). 

\begin{figure}[tbph]
\centering
\includegraphics[width = 2.5in]{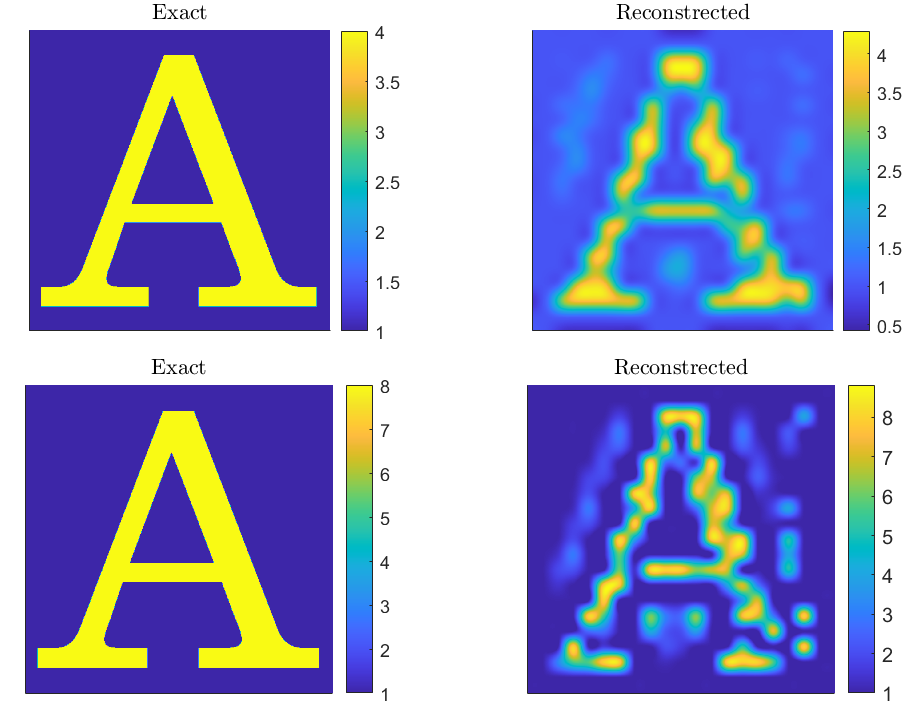}
\caption{Test 2. Exact (left) and reconstructed (right) coefficient $\protect%
\sigma (\mathbf{x})$ with $\protect\sigma_{a}=4$ (first row) and $\protect%
\sigma_{a}=8$ (second row) inside of the letter `A' as in (\protect\ref{8.01}%
). The inclusion/background contrasts in (\protect\ref{8.02}) are
respectively $4:1$ and $8:1$. The reconstructions of both shapes of
inclusions and the inclusion/background contrasts (\protect\ref{8.02}) are
accurate, although the image in the second row is more blurred.}
\label{plot_re_A_4_8}
\end{figure}

\textbf{Test 2.} We test the case when the inclusion in (\ref{8.01}) has the
shape of the letter `$A$' for different values of the parameter $\sigma
_{a}=4,8$ inside of the letter `$A$'. Hence, by (\ref{8.02}) the
inclusion/background contrasts now are respectively $4:1$ and $8:1$.
Computational results are displayed on Figure \ref{plot_re_A_4_8}. One can
observe that shapes of inclusions are imaged accurately. In addition, the
computed inclusion/background contrasts (\ref{8.02}) are accurate. 

\begin{figure}[tbph]
\centering
\includegraphics[width = 3in]{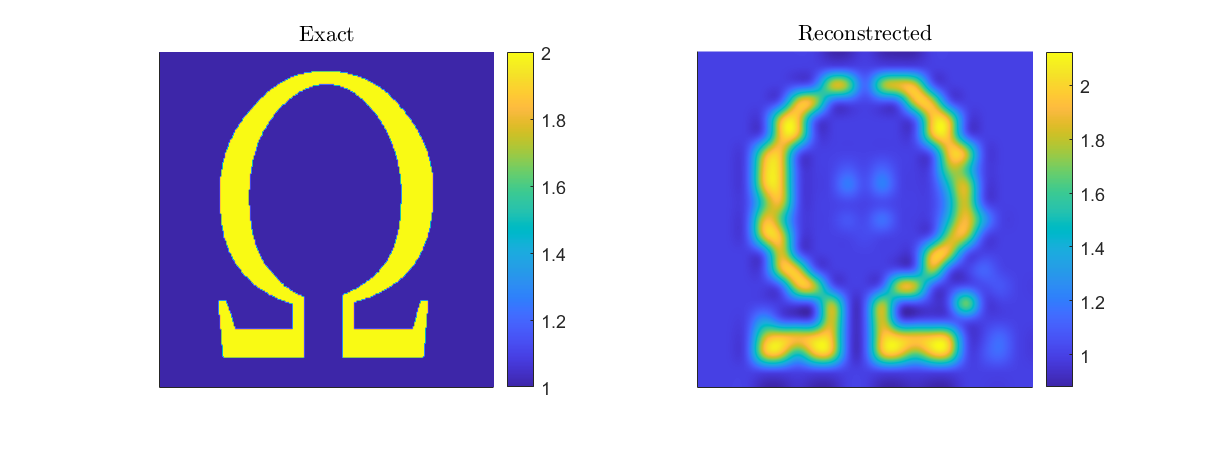}
\caption{Test 3. Exact (left) and reconstructed (right) coefficient $\protect%
\sigma (\mathbf{x})$, where the function $\protect\sigma (\mathbf{x})$ is
given in \eref{8.01} with $\protect\sigma_{a}=2$ inside of the letter `$%
\Omega $'. The reconstructions of both the shape of the inclusions and the
inclusion/background contrast (\protect\ref{8.02}) are accurate.}
\label{plot_re_Omega}
\end{figure}

\textbf{Test 3.} We test the case when the coefficient $\sigma (\mathbf{x})$
in (\ref{8.01}) has the shape of the letter `$\Omega $' with $\sigma _{a}=2$
inside of it. Results are presented on Figure \ref{plot_re_Omega}. We again
observe an accurate reconstruction. 

\textbf{Test 4.} We consider the case when the random noisy is present in
the data in (\ref{8.05}) with $\delta =0.03$, i.e. with the 3\% noise level.
We test the reconstruction for the cases when the inclusion in (\ref{8.01})
has the shape of either the letter `$A$' or the letter `$\Omega $' with $%
\sigma _{a}=2$ inside of them. The results are displayed on Figure \ref%
{plot_A2_Omega2_noisy}. One can observe accurate reconstructions in all four
cases. In particular, the inclusion/background contrasts in (\ref{8.02}) are
reconstructed accurately.

\begin{figure}[tbph]
\centering
\includegraphics[width = 2.5in]{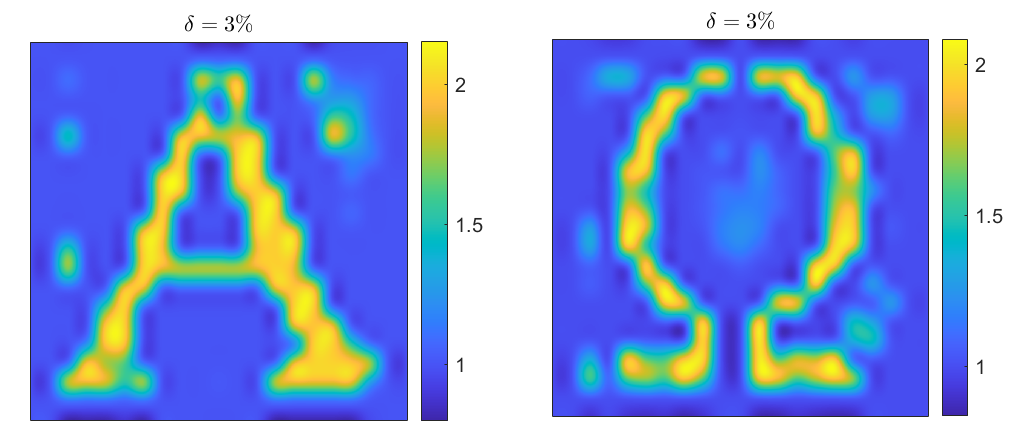}
\caption{Test 4. Reconstructed coefficient $\protect\sigma( \mathbf{x} )$
with the shape of letters `A' (left) or the letter `$\Omega $' (right) with $%
\protect\sigma_{a}=2$ from noisy data \eref{8.05} with $\protect\delta =0.03$%
, i.e. with 3\% noise level. The reconstructions of both shapes of
inclusions and inclusion/background contrasts in \eref{8.02} are accurate. }
\label{plot_A2_Omega2_noisy}
\end{figure}

\newpage

\textbf{Test 5.} In this test, we verify the numerical method for the case,
which is both more complex and more practical one. More precisely, we
consider now the case when the coefficient $\sigma (\mathbf{x})$ in (\ref%
{8.01}) has the shape of a CT scan of an abdomen. We obtain the CT image and
the corresponding segmented image from \cite{EIT_CTScan}. The conductivity
distributions in \cite{EIT_CTScan} range from 2 to 10. However, using a
linear transformation, we obtain the coefficient $\sigma (\mathbf{x})$ in (%
\ref{8.01}) with $\sigma _{a}=2$.

Since this is a completely new set of images, then we need to calibrate our
method again, see item 2 of Remarks 8.1. Hence, we need to select a new
optimal triple $\left( \alpha ,\varepsilon ,\lambda \right) $ of parameters.
We use the same trial and error procedure as the one described above. The
resulting optimal parameters are:%
\begin{equation}
\alpha =0.02,\varepsilon =0.0005,\lambda =2.5.  \label{8.7}
\end{equation}%
The solution of our CIP for this case is presented on Figure \ref%
{plot_re_CT_05_02}, where the left image is exactly the segmented image of 
\cite{EIT_CTScan}. One can observe an accurate reconstruction.

\begin{figure}[htbp]
\centering
\includegraphics[width = 4in]{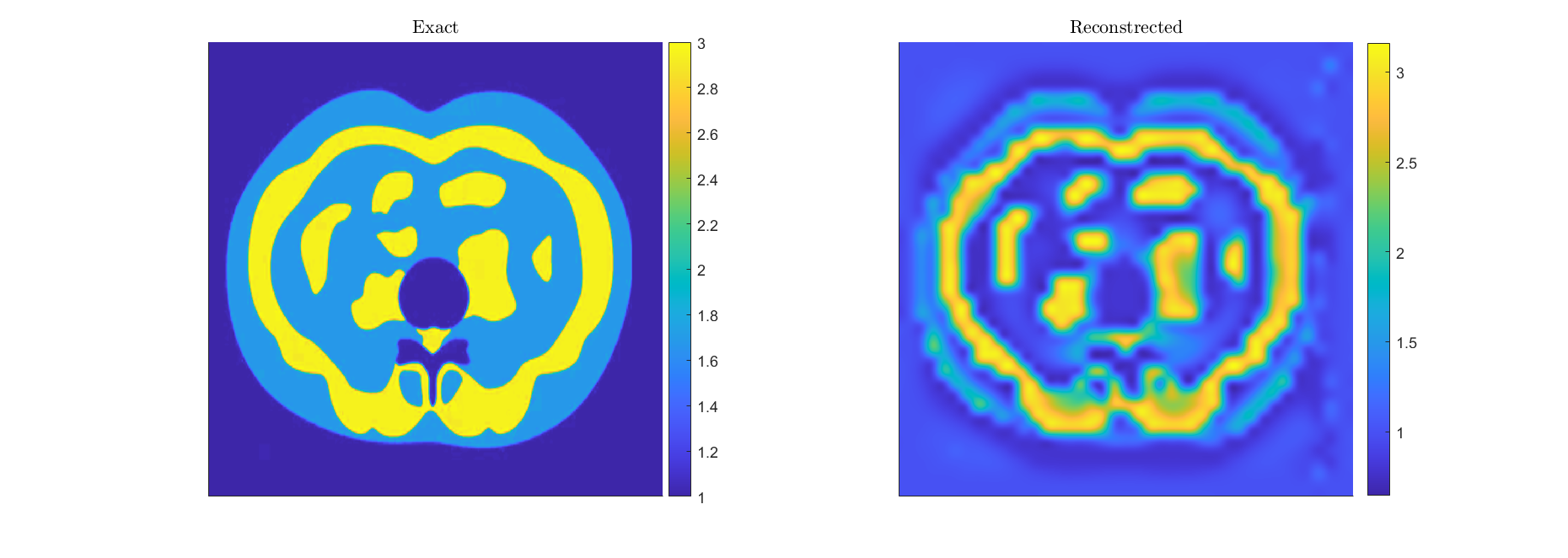}
\caption{Test 5. Abdominal exact segmented image (left) and reconstructed
(right) coefficient $\protect\sigma( \mathbf{x} )$, where the function $%
\protect\sigma( \mathbf{x} )$ is given in \eref{8.01} with $\protect\sigma%
_{a}=3$ with the shape of an abdomen. The reconstructions of both the shape
of the inclusions and the inclusion/background contrast (\protect\ref{8.02})
are approximately good.}
\label{plot_re_CT_06_01}
\end{figure}

Next, we test the optimal triple (\ref{8.7}) of parameters for the segmented
CT scan of \cite{EIT_CTScan}, in which, however, we take $\sigma _{a}=3,$
unlike $\sigma _{a}=2$ of the previous case. The reconstructed image is
displayed on Figure \ref{plot_re_CT_06_01}. The reconstruction is accurate
again.

\begin{figure}[htbp]
\centering
\includegraphics[width = 3in]{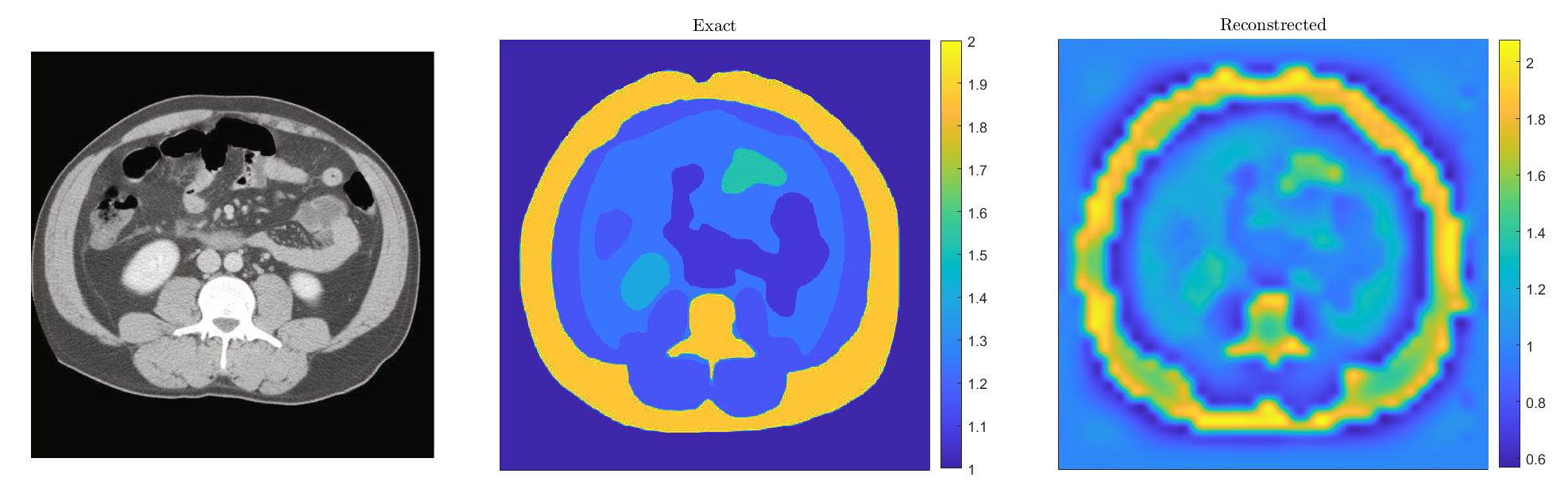}
\caption{Test 5. Abdominal CI image (left), corresponding segmented image
(middle) and reconstructed (right) coefficient $\protect\sigma( \mathbf{x} )$%
, where the function $\protect\sigma( \mathbf{x} )$ is given in \eref{8.01}
with $\protect\sigma_{a}=2$ with the shape of a CT scan of an abdomen. The
reconstruction is rather accurate even for this complicated case with many
colors.}
\label{plot_CT_01_02}
\end{figure}

Finally, we use the set of parameters (\ref{8.7}) to solve our CIP for the
case of a more complicated segmented abdominal image, which we again take
from \cite{EIT_CTScan}. The CT scan is displayed on the left Figure \ref%
{plot_CT_01_02}. The corresponding segmented image is displayed on the
middle Figure \ref{plot_CT_01_02}. One can observe a range of colors here
from dark blue to yellow, which indicates that this case is more complicated
than the ones of two previous images. Our reconstructed image is presented
on the right Figure \ref{plot_CT_01_02}. One can observe that the
reconstruction is rather accurate even for this complicated case with many
colors.

\begin{figure}[htbp]
\centering
\includegraphics[width = 4in]{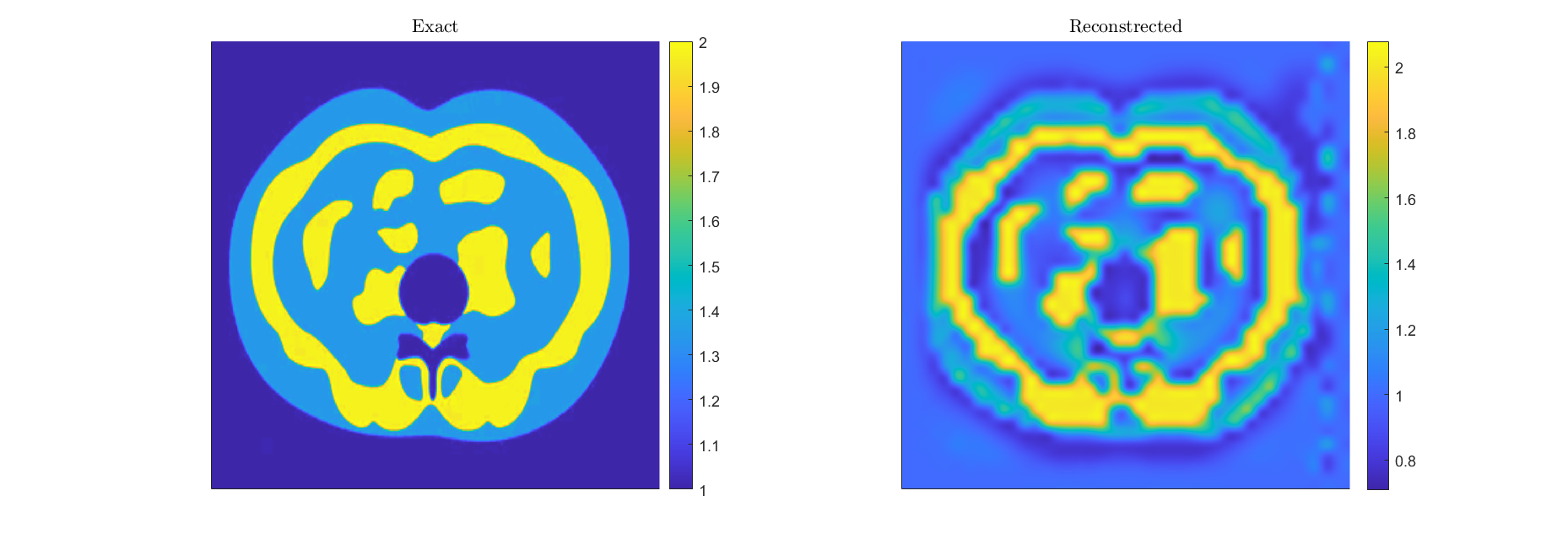}
\caption{Test 5. Abdominal exact segmented image (left) and reconstructed
(right) coefficient $\protect\sigma( \mathbf{x} )$, where the function $%
\protect\sigma( \mathbf{x} )$ is given in \eref{8.01} with $\protect\sigma%
_{a}=2$ with the shape of an abdomen. The reconstructions of both the shape
of the inclusions and the inclusion/background contrast (\protect\ref{8.02})
are accurate.}
\label{plot_re_CT_05_02}
\end{figure}

\section*{Acknowledgement}

The work of MVK was supported by the US National Science Foundation grant
DMS 2436227. The work of Li was partially supported by the Shenzhen Sci-Tech
Fund No. RCJC20200714114556020, Guangdong Basic and Applied Research Fund
No. 2023B1515250005, National Center for Applied Mathematics Shenzhen, and
SUSTech International Center for Mathematics. The work of Yang was partially
supported by NSFC grant 12401558 and Supercomputing Center of Lanzhou
University.

\section*{References}


\end{document}